\documentclass[ejs]{imsart}

\RequirePackage[OT1]{fontenc}
\RequirePackage{amsthm,amsmath}
\RequirePackage[numbers]{natbib}
\RequirePackage[colorlinks,citecolor=blue,urlcolor=blue]{hyperref}
\RequirePackage{etex}
\usepackage{graphicx}

\pubyear{2020}
\volume{0}
\issue{0}
\firstpage{1}
\lastpage{8}
\arxiv{2106.07044}

\usepackage{natbib}%
\usepackage{amsmath,amssymb,mathrsfs}
\usepackage[utf8]{inputenc} %
\usepackage[OT1]{fontenc}    %
\usepackage{hyperref}       %
\usepackage{url}            %
\usepackage{cleveref}
\usepackage{autonum}
\usepackage{booktabs}       %
\usepackage{amsfonts}       %
\usepackage{nicefrac}       %
\usepackage{microtype}      %
\usepackage[dvipsnames,svgnames]{xcolor}         %
\definecolor{darkmidnightblue}{rgb}{0.0, 0.2, 0.4}
\definecolor{darkpowderblue}{rgb}{0.0, 0.2, 0.6}
\definecolor{dukeblue}{rgb}{0.0, 0.0, 0.61}

\hypersetup{
    colorlinks = true,
    citecolor=blue,
    urlcolor=blue,
    breaklinks=true,
    linkcolor = dukeblue,
    linkbordercolor = {white},
}

\usepackage[toc,page]{appendix}
\usepackage{minitoc}

\usepackage{tocloft}

\usepackage{arcs}
\usepackage{graphicx,wasysym}
\usepackage{dsfont,courier}
\usepackage[scaled=.65]{helvet}

\DeclareMathOperator*{\argmin}{arg\,min}

\def\bmu{\boldsymbol\mu}
\def\bY{\boldsymbol Y}
\def\a{\alpha}

\def\z{\zeta}
\def\t{\theta}
\def\T{\Theta}

\def\k{\kappa}

\def\s{\sigma}

\def\O{\Omega}

\def\RR{\mathbb R}

\renewcommand\Im{\operatorname{Im}}

\def\prob{\mathbf P}

\def\bxi{\boldsymbol{\xi}}
\def\btheta{\boldsymbol{\theta}}

\def\bY{\boldsymbol{Y}}

\def\bSigmadiese{\boldsymbol{\Sigma}^\text{\tt\#}}
\def\bX{\boldsymbol{X}}
\def\bXdiese{\boldsymbol{X}^{\text{\tt\#}}}

\def\bthetadiese{\boldsymbol{\theta}^\texttt{\#}}

\def\Xdiese{X^{\text{\tt\#}}}

\def\sdiese{\s^{\text{\tt\#}}}
\def\sdiesesqr{\s^{\text{\tt\#}2}}

\def\xidiese{\xi^\text{\tt\#}}

\def\bsigma{\boldsymbol{\sigma}}

\def\id{\textit{id}}
\def\tdiese{\t^\text{\tt\#}}

\def\btdiese{\boldsymbol{\theta}^\text{\tt\#}}
\def\tdiese{{\theta}^\text{\tt\#}}
\def\bsigmadiese{\boldsymbol{\sigma}^\text{\tt\#}}
\def\bsigmadiesehat{\boldsymbol{\hat\sigma}^\text{\tt\#}}

\def\bQ{\mathbf Q}

\def\kinin{\bar\k_{\textup{in-in}}}
\def\kinout{\bar\k_{\textup{in-out}}}

\newtheorem{theorem}{Theorem}
\newtheorem{lemma}{Lemma}

\newtheorem{remark}{Remark}

\begin{document}

\newenvironment{proofof}[1][Proof]{\begin{trivlist}
\item[\hskip \labelsep {\bfseries #1}]}{\end{trivlist}}

\begin{frontmatter}
\title{Optimal detection of the feature matching map in 
presence of noise and outliers}
\runtitle{Feature matching with outliers}

\begin{aug}
\author{\fnms{Tigran} \snm{Galstyan}\ead[label=e1]{tigran@yerevann.com}}

\address{Russian-Armenian University, YerevaNN\\
9 Charents street,  
0025 Yerevan Armenia\\
\printead{e1}}

\author{\fnms{Arshak} \snm{Minasyan}
\ead[label=e2]{minasyan@yerevann.com}
}

\address{Yerevan State University, YerevaNN\\
9 Charents street,  
0025 Yerevan\\
\printead{e2}}

\author{\fnms{Arnak S.} \snm{Dalalyan}
\ead[label=e3]{arnak.dalalyan@ensae.fr}
\ead[label=u1,url]{www.foo.com}}

\address{CREST, ENSAE, IP Paris\\
5 av.\ Le Chatelier,  
91764 Palaiseau \\
\printead{e3}}

\runauthor{T. Galstyan, A. Minasyan, A. Dalayan}

\end{aug}

\begin{abstract}
We consider the problem of finding the matching map between two sets of $d$-dimensional vectors from noisy observations, where the second set contains outliers.The matching map is then an injection, which can be consistently detected only if the vectors of the second set are well separated. The main result shows that, in the high-dimensional setting, a detection region of unknown injection may be characterized by the sets of vectors for which the inlier-inlier distance is of order at least $d^{1/4}$ and the inlier-outlier distance is of order at least $d^{1/2}$. These rates are achieved using the matching minimizing the sum of logarithms of distances between matched pairs of points. We also prove lower bounds establishing optimality of these rates. Finally, we report the results of numerical experiments on both synthetic and real world data that illustrate our theoretical results and provide further insight into the properties of the algorithms studied in this work.
\end{abstract}


\end{frontmatter}

\section{Introduction}\label{sec:intro}

Finding the best match between two clouds of 
points is a problem encountered in many real
problems. In computer vision, one can look for
correspondences between two sets of local
descriptors extracted from two images. In 
text analysis, one can be interested in matching
vector representations of the words of two
similar texts, potentially in two different 
languages. The goal of the present work is to
gain theoretical understanding of the statistical
limits of the matching problem.  

In the sequel, we use the notation $[n] = \{1,\ldots, 
n\}$ for any integer $n$, and define $\|\cdot\|$ 
as the Euclidean norm in $\mathbb{R}^d$. Assume 
that two independent sequences $\bX = (X_i;
i\in[n])$ and $\bY =(Y_i; i\in[n])$ of independent 
vectors are generated such that $X_i$ and $Y_i$ 
are drawn from the same distribution $P_i$ 
on $\mathbb{R}^d$, for every $i\in[n]$. The 
statistician observes the sequence $\bX$ 
and a shuffled version
$\bXdiese$ of the sequence $\bY$. More precisely, 
$\bXdiese$ is such that $\Xdiese_i = Y_{\pi^*(i)}$
for some unobserved permutation $\pi^*$. The 
goal of matching is to infer the permutation
$\pi^*$ from data $(\bX,\bXdiese)$. In the case
of Gaussian distributions $P_i$, this problem
has been studied in \citep{jmlr_CD13, collier2016minimax}. Clearly, consistent detection
of the matching map $\pi^*$ is impossible if
there are two data generating distributions
$P_i$ and $P_j$ that are very close. In 
\citep{jmlr_CD13, collier2016minimax}, a 
precise quantification of the separation 
between these distributions is given that 
enables consistent detection of $\pi^*$. 
Furthermore, it is shown that the permutation
minimizing the sum of logarithms of pairwise 
distances between the elements of 
$\bX$ and the elements of the shuffled version $\bXdiese$ is an optimal detector of $\pi^*$. 

In this paper, we extend the model studied in \citep{collier2016minimax} to the case when the
set $\bXdiese$ is contaminated by outliers. The number of outliers is supposed to be known and is equal to $m - n$, where $n = |\bX|$ and $m = |\bXdiese|$ are the sizes of considered 
two sequences, however the indices of the 
outliers are unknown. All the distributions 
are assumed in this paper to be
spherical Gaussian, although all the probabilistic
tools used in the proofs have their sub-Gaussian
counterparts. Thus, we consider that two spherical
Gaussian distributions
\footnote{We use the notation $\mathbf I_d$ for the $d\times d$ identity matrix}
$P_1 = \mathcal N_d(\mu_1, \sigma_1^2\mathbf I_d)$ and $P_2 = \mathcal
N_d(\mu_2, \sigma_2^2\mathbf I_d)$ are well
separated if the ``distance to noise ratio''
$\kappa(P_1,P_2) = \|\mu_1-\mu_2\|/\sqrt{
\sigma_1^2+\sigma_2^2}$ is large. Main findings
of \citep{collier2016minimax}, in terms of
smallest separation distance $\bar\kappa = 
\min_{i\neq j}\kappa(P_i,P_j)$ are summarized
in the second columns of \Cref{tab:1}. 
Likewise, the last column of the table 
provides a summary of the contributions of
the present paper in terms of $\kinin = 
\min_{i\neq j} \kappa(P_i,P_j)$ and $\kinout = 
\min_{i,j} \kappa(P_i,Q_j)$, where $Q_1,\ldots,
Q_{m-n}$ are the distributions of the outliers. 

An unexpected finding of this work is that the
``degree'' of heteroscedasticity of the model
has a strong impact on the separation distances
and the detection regions (sets of values of 
$(\kinin,\kinout)$ for which it is possible to
detect the feature map $\pi^*$). This is in sharp
contrast with the outlier-free case, where consistent
detection requires $\bar\kappa$ to be 
at least of order $(d\log n)^{1/4}$ irrespective 
from the behaviour of variances of $P_i$. 
We prove in this work that in the high dimensional
regime $d\ge c\log n$, which is arguably more 
appealing than the low dimensional regime $d\le c\log n$,
the following statements are true:
\begin{itemize}
    \item  If there is no heteroscedasticity, \textit{i.e.}, when all the variances are equal, consistent detection
    of $\pi^*$ is possible if and only if $\bar\kappa 
    = \kinin\wedge\kinout $ 
    is at least of order $(d\log(nm))^{1/4}$. This is
    the same rate as in the outlier-free case.
    \item If the heteroscedasticity is mild, \textit{i.e.}, 
    all the variances are of the same order, the 
    condition
    $\kinin\gtrsim (d\log(nm))^{1/4}$ is the same as in the
    previous item, but the stronger condition $\kinout
    \gtrsim d^{1/2}$ is needed for the inlier-outlier
    separation distance.
    \item Finally, in the general heteroscedastic setting
    both $\kinin$ and $\kinout$ should be at least of 
    order $d^{1/2}$. Furthermore, in all these cases
    consistent detection is performed by the same
    procedure: the Least Sum of Logarithms (LSL). 
\end{itemize}
Note also that the empirical evaluation reported in the present 
paper shows that LSL is attractive not only from the
theoretical but also from the practical point of view.

\begin{table}
    \begin{tabular}{ c|c|c } 
        \toprule
         & Without outliers  & With outliers \\
         &\citep{collier2016minimax}& \textit{current paper} \\
        \toprule
        \shortstack{known $\bsigmadiese$ or \\ 
        all equal $\bsigmadiese$s} & \shortstack{ 
        LSNS is optimal \\ $\bar\kappa \gtrsim (d\log n)^{1/4}$} 
        & \shortstack{ LSNS is optimal [Thm. 1] \\ 
        $\kinin \wedge \kinout\gtrsim (d
        \log(nm))^{1/4}$} \\ 
        \midrule
        \begin{tabular}{c}
        \shortstack{unknown $\bsigmadiese$ \\
        $\nicefrac{\sigma_{\max}}{\sigma_{\min}} \le C$} \\
        \midrule
        \shortstack{unknown $\bsigmadiese$ \\
        arbitrary} 
        \end{tabular}
        & \shortstack{LSL is optimal \\[5pt] $\bar\kappa 
        \gtrsim (d\log n)^{1/4}$} & 
        \begin{tabular}{c}
            \shortstack{LSL is optimal [Thm. 4] \\
            $\kinin \gtrsim (d\log(nm))^{1/4}$  \& 
            $\kinout\gtrsim d^{1/2}$}  \\
            \midrule
            \shortstack{ LSL is optimal 
            [Thm. 2, 3] \\ $\kinin \wedge 
            \kinout \gtrsim d^{1/2}$} 
        \end{tabular}
        \\
        \bottomrule
    \end{tabular}
    \caption{A brief overview of the contributions 
    in the high-dimensional regime $d\ge c\log n$.
    The table provides the condition on the normalized 
    inlier-inlier distance $\kinin$ and inlier-outlier 
    distance $\kinout$, making it possible to consistently
    detect the matching map between two sets of 
    $d$-dimensional vectors. LSL and LSNS refer to 
    least sum of logarithms and least sum of normalized
    squares, respectively.
    }

    \vspace{-19pt}
    \label{tab:1}
\end{table}

\paragraph{Agenda} \Cref{sec:model} describes the 
framework of the vector matching problem and 
introduces the ter\-minology
used throughout this paper. Precise statements 
of the main theoretical results are gathered in \Cref{sec:main}. The prior work is briefly discussed in 
\Cref{sec:prior}. \Cref{sec:numerical} contains numerical experiments carried out both for synthetic and real data. A brief summary and some concluding remarks are presented in \Cref{sec:conc}. Proofs of all theoretical claims are deferred to the supplemental material.  

\section{Problem Formulation}\label{sec:model}

{%
We begin with formalizing the problem of matching two 
sequences of feature vectors $(X_1,\ldots,X_n)$ and
$(\Xdiese_1,\ldots,\Xdiese_m)$ with different sizes $n$ 
and $m$ such that $m \ge n \ge 2$. In what follows, we 
assume that the observed feature vectors are randomly 
generated from the model
\begin{equation}\label{model}
\begin{cases}
X_i = \t_i + \s_i\xi_i \ , \\
\Xdiese_j = \tdiese_j + \sdiese_j\xidiese_j,
\end{cases}\quad i=1,\ldots,n \text{ and } 
j=1,\ldots,m.
\end{equation}
In this model, illustrated in \Cref{fig:ani}, it
is assumed that
\vspace{-5pt}
\begin{itemize}\setlength{\itemsep}{-2pt}
\item $\btheta=(\t_1,\ldots,\t_n)$ and $\btdiese =
(\tdiese_1,\ldots,\tdiese_m)$ are two sequences of 
vectors from $\RR^d$, corresponding to the original 
features, which are unavailable,
\item $\bsigma=(\s_1,\ldots,\s_n)^\top, \bsigmadiese 
= (\sdiese_1,\ldots,\sdiese_m)^\top$ are positive real 
numbers corresponding to the magnitudes of the noise
contaminating each feature, 
\item $\xi_1, \ldots, \xi_n$ and $\xidiese_1, \ldots,
\xidiese_m$ are two independent sequences of i.i.d.\ 
random vectors drawn from the Gaussian distribution 
with zero mean and identity covariance matrix.
\end{itemize}

The simplest special case of \eqref{model}, 
considered in \citep{collier2016minimax}, corresponds 
to the situation where a perfect matching exists
between the two sequences $\btheta$ and $\btdiese$. 
This means that $m = n$ and, for some bijective
mapping $\pi^*:[n]\to[n]$, $\theta_i = \tdiese_{
\pi(i)}$ for all $i\in[n]$. In the general 
case, both $\bX = (X_1, \dots, X_n)$ and 
$\bXdiese = (\Xdiese_1, \dots, \Xdiese_m)$ may 
contain outliers, \textit{i.e.} feature vectors 
that have no corresponding pair. In such a 
situation, it is merely assumed that there exists 
a set $S\subset [n]$ and an  injective mapping 
$\pi^*:S\to [m]$ such that 
\begin{align}\label{pistar}
    \theta_i = \tdiese_{\pi^*(i)} \quad\text 
    {and}\quad \sigma_i = \sdiese_{\pi^*(i)}, 
    \qquad \forall\, 
    i\in S. 
\end{align}
In this case we say that the vectors $\{X_i:i\in[n]\setminus S\}$ and $\{\Xdiese_j:j\in[m]\setminus\pi^*(S)\}$ are 
outliers. The ultimate goal is to detect the
feature matching map $\pi^*$. 

In this work we consider the case when $S = [n]$ and $m>n$. This means 
that only the larger set of feature vectors, namely 
$\bXdiese$, contains outliers. Let us also define the 
set $O_{\pi^*} \triangleq [m] \setminus 
\text{Im}(\pi^*)$,  which contains the indices of 
outliers and satisfies  $|O_{\pi^*}| = m - n$. 
Naturally, the feature vectors contained in $\bX$, 
as well as those vectors from $\bXdiese$ that are 
not outliers, are called {\it{inliers}}.

In this formulation, the data generating distribution is defined by the parameters
$\bthetadiese$,
$\bsigmadiese$ and $\pi^*$. We omit the set of parameters $\btheta$ and $\bsigma$, since they are automatically identified using $\pi^*$, $\bthetadiese$ and $\bsigmadiese$ by the formula $(\theta_i, \sigma_i) = (\tdiese_{\pi^*(i)},\sdiese_{\pi^*(i)})$ for $i\in[n]$. Since our goal is to match the
feature vectors, we focus our attention on the problem of
detecting the parameter $\pi^*$ only, considering
$\bthetadiese$ and $\bsigmadiese$ as nuisance parameters. 
In what follows, we denote by $\prob_{ \bthetadiese,
\bsigmadiese, \pi^*}$ the probability distribution of the
sequence $(X_1,\ldots,X_n,\Xdiese_1,\ldots,\Xdiese_m)$ 
defined by (\ref{model}) under condition (\ref{pistar}) 
with $S= [n]$. 

We are interested in designing estimators that have an expected error smaller than a prescribed level $\alpha$ under the weakest possible conditions on the nuisance parameter $\bthetadiese$ and noise level $\bsigmadiese$. Clearly, the problem of matching becomes more difficult with hardly distinguishable features. To quantify this phenomenon, we introduce the normalized separation distance $\kinin = \kinin(\bthetadiese, \bsigmadiese, \pi^*)$ and the normalized outlier separation distance $\kinout = \kinout(\bthetadiese, \bsigmadiese, \pi^*)$, which measure the minimal distance-to-noise ratio between inliers and the minimal distance-to-noise ratio between inliers and outliers, respectively. The precise 
definitions read as
\begin{align}\label{kappa}
\kinin \triangleq \min_{\substack{ i,j\not\in O_{\pi^*}, 
\\ j \not=i}} \frac{\|\tdiese_i-\tdiese_j\|}{
(\sdiese_i{}^2+\sdiese_j{}^2)^{1/2}},\qquad
\kinout \triangleq \min_{\substack{ i\not\in O_{\pi^*}, 
\\ j \in O_{\pi^*}}} \frac{\|\tdiese_i-\tdiese_j\|}{
(\sdiese_i{}^2 +\sdiese_j{}^2)^{1/2}}. 
\end{align}
Notice that $\kinin$ can be rewritten as 
\begin{align}
    \kinin = \min_{\substack{i, j \in [n] \\ i \neq j}} \frac{\| \theta_i - \theta_j\|}{(\s_i^2 + \s_j^2)^{1/2}}.
\end{align}

\begin{figure}
\centering
\includegraphics[height=0.32\textwidth]{
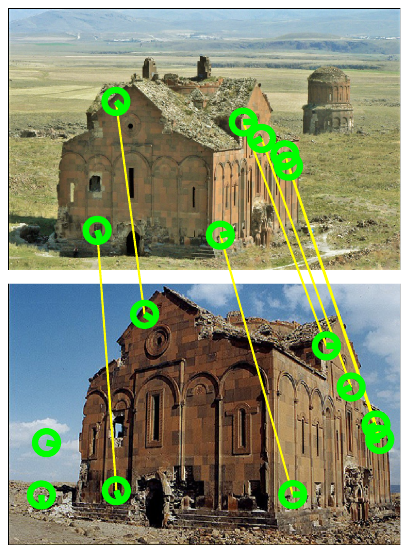}
\hspace{20pt}
\vspace{-11pt}
\includegraphics[height=0.32\textwidth]{
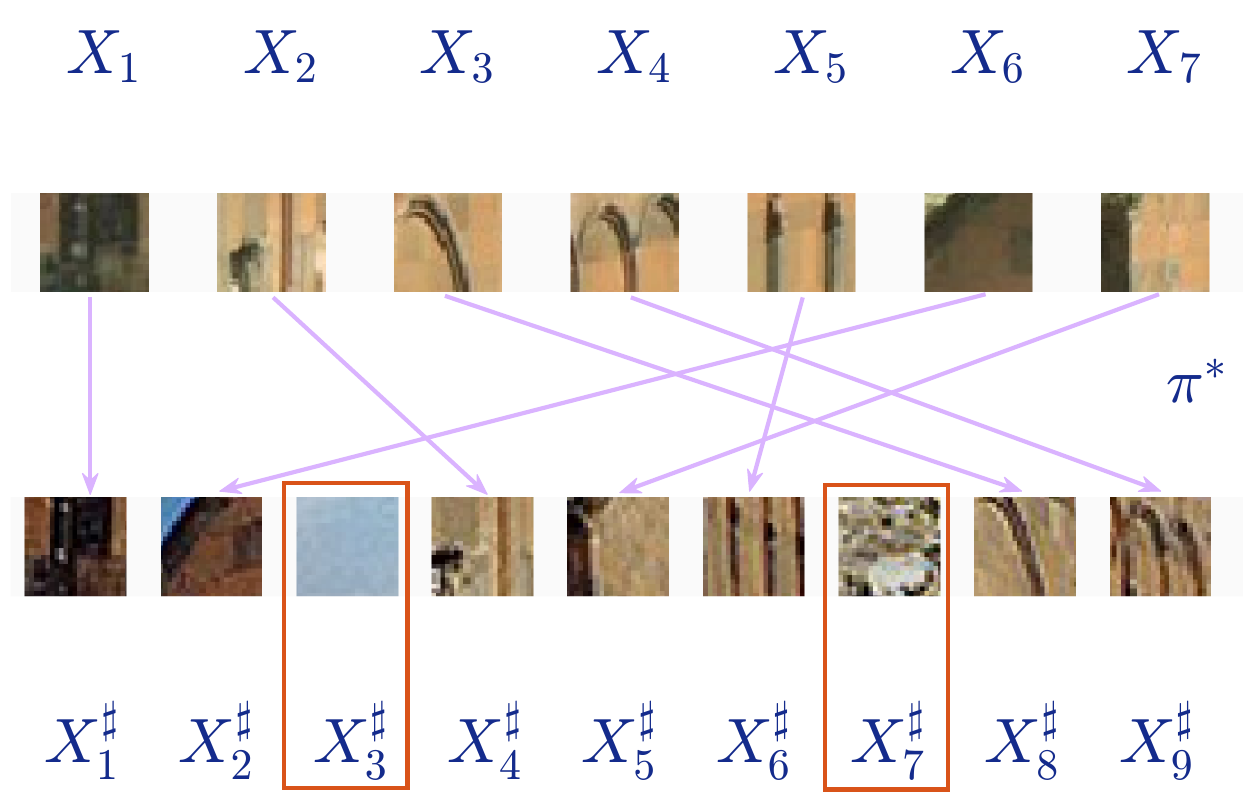}
\caption{Illustration of the considered framework described in \eqref{model}. 
We wish to match a set of 7 patches extracted from
the first image to the 9 patches from the second 
image. The picture on the left shows the locations 
of patches as well as the true matching map $\pi^*$ 
(the yellow lines).
}
\vspace{-13pt}
\label{fig:ani}
\end{figure}

Clearly, if $\kinin=0$ or, $\kinout=0$, there are
two identical feature vectors in $\bXdiese$. In such
a situation, assuming $\sigma_i$'s are all equal, 
the parameter $\pi^*$ is nonidentifiable, in the 
sense that there exist two different permutations
$\pi^*_1$ and $\pi^*_2$ such that the distributions
$\prob_{\bthetadiese,\bsigmadiese,\pi^*_1}$ and
$\prob_{\bthetadiese,\bsigmadiese,\pi^*_2}$ 
coincide.  Therefore, to ensure the existence of
consistent detectors of $\pi^*$ it is necessary 
to impose the conditions $\kinin>0$ and 
$\kinout>0$. Moreover, good procedures are those
consistently detecting $\pi^*$ even if either
$\kinin$ or $\kinout$ are small. We are interested
here in finding the detection boundary in terms 
of the order of magnitude of $(\kinin, \kinout)
$. More precisely, for any given $\alpha\in(0,1)$ 
we wish to find a region $\mathcal R_{n,m,d}^\alpha$
in $\mathbb R^2$ such that:
\setlength{\leftmargini}{2em}
\begin{itemize}
    \item There is an estimator $\hat\pi_{n,m}$ of 
    $\pi^*$ satisfying 
        $\prob_{\bthetadiese, \bsigmadiese, \pi^*}(\hat\pi 
        \neq \pi^*) \le \alpha$
    for every $(\bthetadiese,\bsigmadiese,\pi^*)$ 
    lying in the detection region, \textit{i.e.}, for 
    which $(\kinin, \bar\k_{ 
    \textup{in-out}}) \in \mathcal R_{n,m,d}^\alpha$.
    \item There is a constant $C<1$ such that for any
    estimator $\bar\pi_{n,m}$ of $\pi^*$, we can find
    a parameter value $(\bthetadiese,\bsigmadiese,\pi^*)$
    in the region $\{(\bthetadiese,\bsigmadiese,\pi^*) :
    (\kinin, \kinout) 
    \in C\mathcal R_{n,m,d}^\alpha\}$
    such that $\bar\pi$ fails to detect $\pi^*$ with a
    probability larger than $\alpha$. 
\end{itemize}
Let us make two remarks. First, note that in the 
outlier-free case considered 
in \citep{collier2016minimax}, $\kinout$
is meaningless and, therefore, the detection region  is one-dimensional $\mathcal R_{n,m,d}^\alpha$. Thus, it 
is necessarily a half-line and is proven to be of the
form $\kinin\ge C(\log \nicefrac{ 
n}{\alpha})^{1/2} \vee (d\log \nicefrac{n}{\alpha} )^{1/4}$
for some universal constant $C$. Second, the aforementioned definition of the detection region $\mathcal R_{n,m,d}^\alpha$ does not guarantee
its uniqueness (even up to a scaling by a universal 
constant). This is in contrast with the 
outlier-free case. To overcome this difficulty, we 
look for $\mathcal R_{n,m,d}^\alpha$ of
the form $[t_{\textup{in-in}},+\infty)\times [t_{\textup{in-out}},+\infty)$ with the smallest possible
threshold $t_{\textup{in-out}}$ for the normalized inlier-outlier distance $\kinout$.

\def\bXi{\boldsymbol{\Xi}}
\section{Main theoretical results}\label{sec:main}

In this section, we have collected the main theoretical findings of the paper.  When the noise is homoscedastic,
\textit{i.e.},  when all $\bsigma$'s are equal, 
the results obtained by \cite{collier2016minimax} in the outlier-free setting can be easily extended to the setting
with outliers. Therefore, in the present paper, we focus 
on the heteroscedastic case. For the sake of clarity of
exposition, we will present the results in the case of known
variances $\bsigma,\bsigmadiese$ prior to investigating
the more interesting case of unknown variances.  

The detection regions we study below are based on the 
 maximum profile likelihood estimator. The model presented in \eqref{model} has the parameter $\bXi = (\btdiese, \bsigmadiese, \pi)$, while the observations are the sequences
of feature vectors $\bX$ and $\bXdiese$. The negative
log-likelihood of this model is given by
\begin{align}
    \ell_n(\boldsymbol{\Xi}; \{\bX, \bXdiese\}) = &\sum_{i=1}^n \bigg(\frac{\| X_i - \tdiese_{\pi(i)}\|_2^2}{2\sdiesesqr_{\pi(i)}} + \frac 1 2 \log(\sdiesesqr_{\pi(i)})\bigg)\\ &+  \sum_{j=1}^m \bigg(\frac{\| \Xdiese_j - \tdiese_j\|_2^2}{2\sdiesesqr_j} + \frac 1 2 \log(\sdiesesqr_j)\bigg). 
\end{align}
The profile negative log-likelihood is then defined as 
the minimum
with respect to $(\bthetadiese,\bsigmadiese)$ of 
the log-likelihood $\ell_n(\boldsymbol{\Xi}; \{\bX, \bXdiese\})$. 

\subsection{Warming up: known variances $\bsigma,\bsigmadiese$}

One can check
that the minimization with respect to $\bthetadiese$ leads
to the variance-dependent cost function
\begin{align}\label{liklihood-sigma}
    \ell_n(\pi,\bsigmadiese; \{\bX, \bXdiese\}) =  \sum_{i=1}^n \frac{\|X_{i}-\Xdiese_{\pi(i)}\|^2}{\s_{i}^2+\sdiesesqr_{\pi(i)}} + \sum_{i=1}^n \frac{1}{2}\log(\sdiesesqr_{\pi(i)}) + 
    \sum_{j=1}^m \frac{1}{2} \log(\sdiesesqr_j).
\end{align}
When $m=n$ and there is no outlier, the last two sums of the
last display do not depend on $\pi$ and, therefore, the
maximum profile likelihood estimator of $\pi^*$ is obtained
by the Least Sum of Normalized Squares (LSNS) criterion
\begin{align}\label{LSNS}
    \hat\pi^{\textrm{LSNS}}_{n,m}\in \argmin_{\pi:[n]\to[m]}
    \sum_{i=1}^n \frac{\|X_{i}-\Xdiese_{\pi(i)}\|^2}{\s_{i}^2 
    +\sdiesesqr_{\pi(i)}},
\end{align}
where the minimum is over all injective mappings 
$\pi:[n]\to[m]$. This, and the other estimators defined
in this work, can be efficiently computed using suitable 
versions of the Hungarian algorithm \citep{Kuhn,Kuhn2010,
Munkres}. As shows the next theorem, it 
turns out that even when $m>n$, the estimator
$\hat\pi^{\textrm{LSNS}}_{ n,m } $  defined above 
leads to an optimal detection region. 

\begin{theorem}[Upper bound for LSNS]\label{thm:upperLSNS}
Let $\alpha \in (0,1)$ and condition \eqref{pistar} 
be fulfilled. If the separation distances $\kinin$ 
and $\kinout$ corresponding to $(\bthetadiese, 
\bsigmadiese,\pi^*)$ and defined in \eqref{kappa}
satisfy the condition
\begin{align}\label{thresh:1}
    \min\{\kinin, \kinout\} \ge 4 \Big\{ \big(
    d\log(\nicefrac{4nm}{\alpha})\big)^{1/4} \vee
    \big( 2\log(\nicefrac{8nm}{ \alpha})\big)^{1/2} 
    \Big\}
\end{align}
then the LSNS estimator defined in \eqref{LSNS} detects
the true matching map $\pi^*$ with probability at least
$1-\alpha$, that is
\begin{align}
    \prob_{\bthetadiese, \bsigmadiese, \pi^*}(\hat\pi^{
    \textup{LSNS}}_{n,m} = \pi^*) \ge 1-\alpha.
\end{align}
\vspace{-15pt}
\end{theorem}

The similarity---both its statement and its proof--- 
of this result to its counterpart 
in the outlier-free setting might suggest that 
the presence of outliers does not make the problem 
any harder from a statistical point of view. 
However, this is not true in the more
appealing setting of unknown variances.  

\begin{remark}
All the results of this subsection apply in the 
homoscedastic case---when $\sigma_i=\sdiese_j$ 
for all $(i,j)\in [n]\times[m]$---with unknown 
noise level. Indeed, in this case, the LSNS 
estimator coincides with the minimizer of the 
sum of squared errors and, therefore, does not 
depend on the noise levels. 
\end{remark}

\begin{remark}\label{rem:2}
\Cref{thm:upperLSNS} can be readily extended to 
the case where the noise vectors $X_i-\theta_i$ 
and $\Xdiese_j-\tdiese_j$ are Gaussian with zero 
mean and general covariance matrices, denoted 
respectively by $\boldsymbol\Sigma_i$ and 
$\bSigmadiese_j$. Then, the LSNS estimator should 
be defined as the minimizer of the sum over $i$ 
of the terms $(X_i-\Xdiese_{\pi(i)})^\top 
(\boldsymbol\Sigma_i + \bSigmadiese_{\pi(i)})^{-1}
(X_i-\Xdiese_{\pi(i)})$. Hence, redefining 
\begin{align}\label{kappa_bis}
    \kinin \triangleq 
        \min_{i,j\not\in O_{\pi^*},j \not=i}  
        \kappa_{ij},\qquad 
    \kinout \triangleq 
        \min_{i\not\in O_{\pi^*},j \in O_{\pi^*}} 
        \kappa_{ij},
\end{align}
where $\kappa_{ij}^2 = (\tdiese_i-\tdiese_j)^\top 
(\bSigmadiese_i + \bSigmadiese_j)^{-1} (\tdiese_i 
-\tdiese_j)$, one gets exactly the result as the
one stated in \Cref{thm:upperLSNS}. 
\end{remark}

\begin{remark}
The minimax setting considered in the present work is 
largely inspired by the corresponding setting in the 
problem of statistical hypotheses testing 
\cite{Ingster, JudNem,Yuting2,Yuting,carpentier2019minimax}. It should be noted that in 
many hypothesis testing problems, one can further 
the results on minimax rates of separation by obtaining
sharp constants \citep{Ermakov,lepski2000asymptotically,ComDal2}. It would
be interesting to investigate whether it is possible or
not to obtain sharp constants in the setting of this
paper. To the best of our knowledge, this question is
open for other instances of multiple hypotheses testing 
as well.
\end{remark}

\subsection{Detection of $\pi^*$ for unknown and 
arbitrary variances $\bsigma,\bsigmadiese$}

The LSNS procedure analyzed in \Cref{thm:upperLSNS} 
exploits the values of known noise variances to normalize 
the squares of distances between vectors $X_i$ and 
$\Xdiese_{\pi(i)}$. Therefore, LSNS is inapplicable 
in the case of unknown noise variances, unless the values of all $\sigma_i$ and $\sdiese_i$ are equal. 
Instead, we consider the Least Sum of 
Logarithms (LSL) estimator
\begin{align}\label{LSL}
    \hat\pi^{\textup{LSL}}_{n,m} \triangleq \argmin_{\pi:[n]\to [m]}\sum_{i=1}^n {\log \|X_{i}-\Xdiese_{\pi(i)}\|^2},
\end{align}
where the minimum is over all injective maps $\pi
: [n] \to [m]$. This estimator can be seen as the
minimizer of a criterion defined as the minimum of
the cost function from \eqref{liklihood-sigma} with
respect to $\bsigmadiese$ under the constraint
$\min_{j\not\in\Im(\pi)} \sdiese_{j} \ge 
\sigma_{\min}$, for some fixed (but unknown) 
constant  $\sigma_{\min}>0$. 

To provide a quick overview of what follows, let us
stick in the remaining of this paragraph to the 
case $\log(nm) = O(d)$ so that the right hand side 
of \eqref{thresh:1} is of order $\big( d\log(nm) 
\big){}^{1/4}$. Recall that in the outlier-free 
case, the LSL estimator has been shown to perform 
as well as the LSNS while having the advantage of 
not requiring the knowledge of variances
$\bsigmadiese$ \citep{collier2016minimax}. 
Somewhat unexpectedly, the situation is significantly
different in the presence of outliers.  Indeed, the
best we managed to prove in the presence of 
outliers is that the detection of the matching 
map by LSL is possible whenever $\min\{ \bar\k_{ 
\textup{in-in}}, \kinout\}\ge
C\sqrt{d}$ for some sufficiently large constant 
$C$. The precise statement being given in the next
theorem, let us mention right away that the
discrepancy between this rate $\sqrt{d}$ and the 
rate $\big( d\log(nm) \big){}^{1/4}$ in
\eqref{thresh:1} is due to the inherent hardness 
of the setting and not merely an artefact of the
proof. This will be made clear below.

\begin{theorem}[Upper bound for LSL]
\label{thm:upper-bound-LSL}
Let $\alpha \in (0,1/2)$ and condition 
\eqref{pistar} be fulfilled. If the separation 
distances $\kinin$ and $\kinout$ corresponding 
to $(\bthetadiese,\bsigmadiese,\pi^*)$ and 
defined by \eqref{kappa} satisfy
\begin{align}\label{detect}
    \min\{\kinin, \kinout\} \ge \sqrt{2d} + 4 
    \Big\{\big(2d\log (\nicefrac{4nm}{ \alpha}) 
    \big)^{1/4} \vee \big(3 \log(\nicefrac{8nm}{
    \alpha}\big)^{1/2}\Big\}
\end{align}
then the LSL estimator \eqref{LSL} detects
the matching map $\pi^*$ with probability at least
$1-\alpha$, that is
\begin{align}
    \prob_{\bthetadiese, \bsigmadiese, \pi^*}
    (\hat\pi^{\textup{LSL}}_{n,m} = \pi^*) 
    \ge 1 - \alpha.
\end{align}
\vspace{-15pt}
\end{theorem}

This result is disappointing since it requires
the distance between different feature vectors to
be larger than $\sqrt{2d}$ in order to be able
to consistently detect the matching map $\pi^*$. 
As we show below, without any further condition 
(for instance, on the noise variances), this 
rate cannot be improved. Furthermore, the
rate $\sqrt{d}$ is optimal not only for LSL but 
also for the larger class of so called
\textit{distance based $M$-estimators.}

We will say that an estimator $\hat\pi_n$ of $\pi^*$
is a distance based $M$-estimator, if for a sequence of non-decreasing functions $\rho_i:\mathbb{R}_+ 
\to \mathbb{R}$, $i=1,\ldots,n$, the following is correct
\begin{align}
    \hat\pi_n \in \argmin_{\pi :[n] \to [m]} \sum_{i=1}^n \rho_i\big(\| X_i - \Xdiese_{\pi(i)}\|\big),
\end{align} 
where the minimum is over all injective mappings
$\pi:[n]\to[m]$. We denote by $\mathcal M$ the set 
of all distance based $M$-estimators.  
We show that there is indeed a setup where $\kinin \wedge \bar\k_{ \textup{in-out}}$ is as large as $0.2\sqrt{d}$ but
any estimator from $\mathcal M$ fails to detect
$\pi^*$ with probability at least $1/4$.  The next theorem formalizes the described result. 
\begin{theorem}[Lower bound over $\mathcal{M}$]\label{thm:lower-bound}
Assume that $m > n \ge 4$ and $d \ge 422\log(4n)$. There exists a triplet $(\bsigmadiese, \btdiese, \pi^*)$ such that  $\kinin = \kinout = \sqrt{d/20}$ and 
    \begin{align}\label{lower:1}
        \inf_{\hat\pi \in \mathcal{M}} 
        \prob_{\bthetadiese, \bsigmadiese, \pi^*}(\hat\pi 
        \neq \pi^*) > \nicefrac 1 4.
    \end{align}
    \vspace{-15pt}
\end{theorem}
The proof of this theorem, postponed to the appendix,
is constructive. This means that we exhibit a triplet 
$(\bsigmadiese, \btdiese, \pi^*)$ satisfying 
\eqref{lower:1}. Careful inspection shows that in the 
case $d = O(\log(nm))$ the same triplet satisfies 
$\kinin \wedge \bar\k_{ \textup{in-out}} \asymp 
\sqrt{\log(nm)}$ and \eqref{lower:1} is still true. 
This implies that the order of magnitude of the right
hand side of \eqref{detect} is optimal both in the
high-dimensional regime $d\ge 422\log(4n) $ and in the
low-dimensional regime $d< 422\log(4n)$. This shows the optimality of LSL among all estimators 
from $\mathcal{M}$. Note that the estimator $\hat\pi^{
\textup{LSNS}}_{n,m}$ does not belong to the family of 
\textit{distance based M-estimators}. Furthermore, in the
low dimensional regime $d=O(\log(nm))$, the separation 
rate of the LSL, $\sqrt{\log(nm)}$, is the same as that of the LSNS.

The next theorem extends the result of \Cref{thm:lower-bound} establishing the lower bound over all injective mappings $\pi: [n] \to [m]$, hence implying the optimality of the rate presented in \Cref{thm:upperLSNS}. We show that even if $\kinin$ and $\kinout$ are of order $(d\log(nm))^{1/4} \vee (\log(nm))^{1/2}$ then there are indeed  scenarios in which any estimator $\hat\pi$ fails to detect 
$\pi^*$ with probability at least $1/3$.

\begin{theorem}[General lower bound]\label{thm:general-lb}
Denote $\kappa = \min\{\kinin, \kinout\}$. Assume that $m > n \ge 5$ and $d \ge 16\log(nm)$. Then, there exists a triplet $(\bsigmadiese, \btdiese, \pi^*)$ such that $6\kappa \ge (d \log(nm))^{1/4}$ and
\begin{align}
    \inf_{\hat\pi} \prob_{\btdiese, \bsigmadiese, \pi^*} (\hat\pi \neq \pi^*) > \nicefrac{1}{3},
\end{align}
where the infimum is taken over all injective matching maps $\pi:[n] \to [m]$. 
\end{theorem}

In the next section we show that under some mild
conditions on $\bsigmadiese$ it is indeed possible to
obtain different rates for $\kinin$
and $\bar\k_{ \textup{in-out}}$, namely we show that
if $\kinin \gtrsim d^{1/4}$ and
$\kinout \gtrsim d^{1/2}$ then the
LSL estimator detects correct matching with high
probability.

\subsection{Detection of $\pi^*$ for unknown and mildly 
varying variances %
}

The results of the last two theorems are disappointing, 
since they indicate that the features should be very 
different from one another for detection of the matching 
map to be possible. An interesting finding, presented 
below, is that strong constraint can be significantly 
alleviated in the context of mild heteroscedasticity. 
By mild heteroscedasticity we understand here the 
situation in which all variances $\sdiese_i$ are of 
the same order of magnitude.  

\begin{theorem}[Upper bound under mild 
heteroscedasticity]\label{thm:mild}
Let $r_\sigma = \max_{i,j}(\sdiese_i/\sdiese_j)$. 
If the separation distances $\kinin$
and $\kinout$ defined in \eqref{kappa} satisfy
\begin{align}
    \kinin &\ge 
    2\big(4d\log(\nicefrac{4nm}{\alpha}) 
    \big)^{1/4} + 2\big(2\log(\nicefrac{4nm}{
    \alpha})\big)^{1/2}\\
    \kinout 
        &\ge\  \sqrt{2(r_\sigma-1)d} + 2\big(4
        r_\sigma^2 d\log(\nicefrac{4nm}{\alpha}
        ) \big)^{1/4} + 2\big(2 r_\sigma \log( 
        \nicefrac{4nm}{\alpha}) \big)^{1/2} ,
\end{align}
then the LSL estimator \eqref{LSL} 
detects the matching map $\pi^*$ with probability 
at least
$1-\alpha$, that is
\begin{align}
    \prob_{\bthetadiese, \bsigmadiese, \pi^*} 
    (\hat\pi^{\textup{LSL}}_{n,m} = \pi^*) \ge 
    1 - \alpha.
\end{align}
\vspace{-15pt}
\end{theorem}
Note that a lower bound similar to that of
\Cref{thm:lower-bound} can be proved in the case
of mild heterescodestacity as well, showing that 
there is an example for which $\kinin$ is of order
$d^{1/4}$, $\kinout$ is of order $d^{1/2}$ and
any estimator from $\mathcal M$ fails to detect
$\pi^*$ with probability at least $1/4$. 

We complete this section by summarizing the 
joint contribution of 
\Cref{thm:upperLSNS,thm:upper-bound-LSL,thm:lower-bound,thm:mild}. 
To simplify this discussion, we consider two 
cases: high-dimensional case refers to $d\ge \log(4nm/
\alpha)$ (presented in \Cref{tab:1}) and low-dimensional case refers to 
the condition $d< \log(4nm/\alpha)$. In the 
high dimensional setting with arbitrary noise 
variances, the detection region for the LSL
estimator is given by $\{\kinin\wedge \kinout
\ge 15 \sqrt{d}\}$, which is much worse than 
the detection region for LSNS, $\{\kinin
\wedge \kinout\ge 8 (d\log(4nm/\alpha))^{1/4} 
\}$, obtained in the known-variance scenario.  
Somewhat surprisingly, in such a setting, 
even a strong assumption on the outliers, such
as requiring them to be at least at a
distance $0.2\sqrt{d}$ of the inliers, is not
enough for relaxing the assumption on the 
inlier-inlier separation distance. Finally, on a positive
note, in the intermediate case of mildly 
varying variances, the detection region for
the LSL estimator is of the form 
$\{\kinin\ge 7 (d\log(4nm/\alpha))^{1/4}; \ 
\kinout \ge 10 \sqrt{d}\}$. This means that
if the outliers are at a distance $\Omega(
\sqrt{d})$ of the inliers, then the LSL 
recovers the true matching under the same
condition on $\kinin$ as in the outlier-free
setting.
}

\begin{figure}%
\begin{center}
\includegraphics[width=0.65\textwidth]{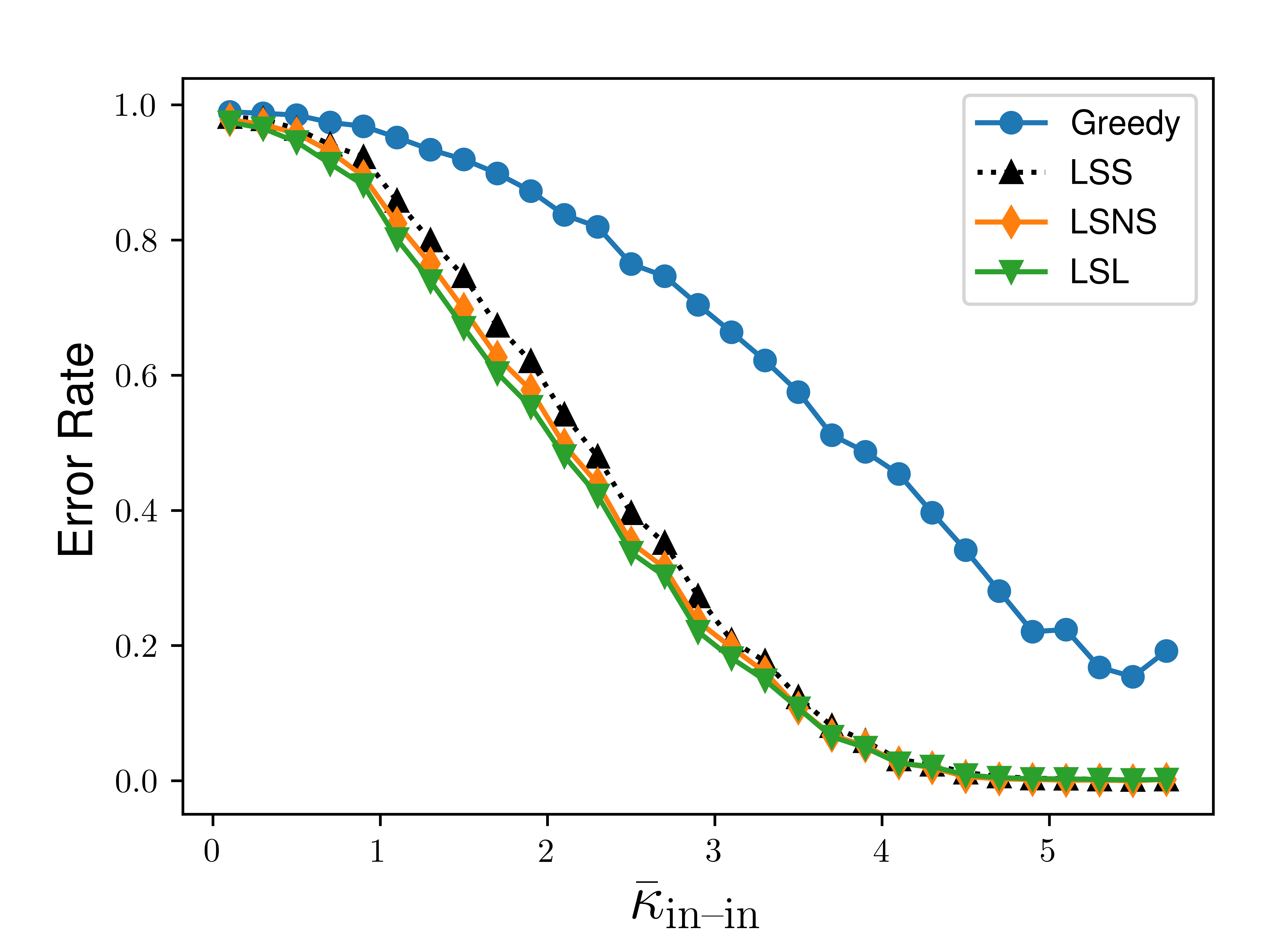}
\vspace{-10pt}
\caption{The performance of the methods (Greedy, LSS, 
LSL, LSNS) in the setup described in Exp.~1. 
Curves represent the error rates (percentage of repetitions
in which the estimated matching map differs from the 
true one) as a function of separation distances. The 
 picture illustrates that LSS, LSL and LSNS require 
much lower value of $\kinin$ in order to find the 
correct mapping.  }
\label{fig:exp1}
\end{center}
\vspace{-20pt}
\end{figure}

\section{Other related work}\label{sec:prior}

Measuring the quality of the various statistical 
procedures of decision making by their minimal 
separation rates became the standard in hypotheses
testing, see the seminal papers \citep{Burnashev,
Ingster82} and the monographs \citep{Ingster,JudNem}. 
Currently this approach is widely adopted in machine
learning literature \citep{Xing20,Wolfer,Balchard,
RamdasISW16,Wei,Collier12a}. Beyond the classical 
setting of two hypotheses, it can also be applied 
to multiple hypotheses testing framework, for instance,
variable selection \citep{Ndaoud,azais2020multiple,
ComDal1} or the matching problem considered here.

On the other hand, feature matching is a well 
studied problem in computer vision. In recent 
years, a great deal of attention was devoted to 
the acceleration of greedy matching algorithms, 
based on approximate and fast methods of finding 
nearest neighbors (e.g. 
\cite{jiang2016, wang2018, wang2011, harwood2016,malkov2020}). 
Another direction that helps to improve feature 
matching problem is using alternative local 
descriptors \citep{orb2011,piifd2010, brief2010} 
for given keypoints. Naturally, the question of 
how to chose keypoints arises, which is addressed, 
for instance, in \citep{bai2020d3feat, Tian_2020_ACCV}. 
For more complete overview of the field we refer 
to \citep{ma2021} and references therein. 

Finally, permutation estimation and related problems 
have been recently investigated in different contexts 
such as statistical seriation 
\citep{flammarion2019optimal}, noisy sorting 
\citep{mao2018minimax}, regression with shuffled data
\citep{pananjady2017linear,slawski2019linear}, 
isotonic regression and matrices 
\citep{Mao18,pananjady2020isotonic,ma2020optimal},
crowd labeling \citep{ShahBW16a}, and recovery of 
general discrete structure \citep{gao2019iterative}.

\section{Numerical results}\label{sec:numerical}
In this section, we report the results of some 
numerical experiments carried out on simulated and real 
data. We applied aforementioned methods LSNS and LSL 
and computed different measures of their performance. 
To get a more complete picture, we included in this
study the Least Sum of Squeres (LSS) estimator 
and the greedy estimator. LSS is an unnormalized 
version of LSNS, given by
\begin{align}\label{LSS}
    \hat\pi^{\textrm{LSS}}_{n,m}\in \argmin_{\pi:[n]\to[m]}
    \sum_{i=1}^n \|X_{i}-\Xdiese_{\pi(i)}\|^2.
\end{align}
It coincides with LSNS in the case of homoscedastaic 
noise. The greedy estimator is obtained by sequentially 
matching each vector from $\bX$ to the nearest vector
from $\bXdiese$. Experiments were implemented using 
python or matlab. For solving linear sum assignment 
problems such as \eqref{LSL} or \eqref{LSS}, the 
generalized and
improved versions of the Hungarian algorithm
were used \citep{Kuhn,Kuhn2010,Munkres,Duff}, implemented in SciPy library \citep{2020SciPy-NMeth}.  The goal of the first
two experiments is to illustrate our theoretical 
findings on synthetic data sets. The third 
experiment aims to highlight that the methods studied in this work have some additional
attractive features that would be interesting 
to investigate in the future.

\paragraph{Experiment 1: Synthetic data with random features}
We first randomly generated $\pi^*$,  
$\bthetadiese$ and $\bsigmadiese$ as follows. We 
randomly sampled from uniform distribution on $[0, 2]$
independent variables $\tau_{ij}$, $i\in[m]$, $j\in[d]$.
Then, $(\tdiese_i)_j$ are independently sampled from 
the Gaussian distribution with 0 mean and variance
$\tau_{ij}$. Additionally, for every $\tdiese_i \in 
\bthetadiese$ such that $i \notin J_{\pi^*}$ 
($\tdiese_i$ is an outlier), we incremented every
coordinate of $\tdiese_i$ by $i$. Entries of 
$\bsigmadiese$ were sampled from the uniform 
distribution over $[0.5, 2]$. Sequences $\bX$ and 
$\bXdiese$ were generated according to \Cref{sec:model} 
with $\pi^*(i) = i$ for $i\in[n]$. We applied to this 
data the following matching algorithms: Greedy, LSS, LSNS and LSL.  

We chose $n = 100$, $m = 130$ and $d = 50$, and 
generated $N=50$ datasets according to the 
foregoing process. For each dataset, we computed
the 0-1 error of the considered estimators and the
values of $(\kinin,\kinout)$.  We plotted in  \Cref{fig:exp1} 
the error averaged over all datasets with a given 
value of $\kinin$. We see that the error 
decreases fast with $\kinin$, corroborating our 
theoretical results. 

\paragraph{Experiment 2: Synthetic data with deterministic features}
The second experiment is conducted on data generated
by features $\bthetadiese$ and variances $\bsigmadiese$
inspired by the example constructed in the proof of
\Cref{thm:lower-bound}. More precisely, for some 
real numbers $a$ and $b$ representing, respectively,
the scale of inlier-inlier distance $\kinin$ and 
inlier-outlier distance $\kinout$, we set
$\tdiese_k = [ka, 0,\ldots,0]^\top$ for 
$k\in[n]$ and $\tdiese_{n+k} = [na + kb, 0,\ldots,
0]^\top$. We also used decreasing variances
$\sdiese_k = 1/k^{3/2}$ for $k\in[m]$ and true identity
mapping $\pi^*(k)=k$ for $k\in [n]$. We chose 
$n=100$, $m=120$ and dimension $d$ varying in
the set $\{10,20,40\}$. For each pair of values
$(a,b)$ in a suitably chosen grid, we repeated
$n_{\textup{rep}} = 400$ times the experiment
that consisted in generating data according to 
\eqref{model} and computing estimators $\hat\pi_{n,m
}^{\textup{LSS}}$ and $\hat\pi_{n,m}^{\textup{LSL}}$
defined respectively by \eqref{LSS} and \eqref{LSL}. 
We then computed, for each pair $(a,b)$ and for each
estimator LSS and LSL, the percentage of successful
detection among $n_{\textup{rep}}$ repetitions. 

\begin{figure}%
    \centering
    \includegraphics[scale=0.8, width=0.85\textwidth]{
    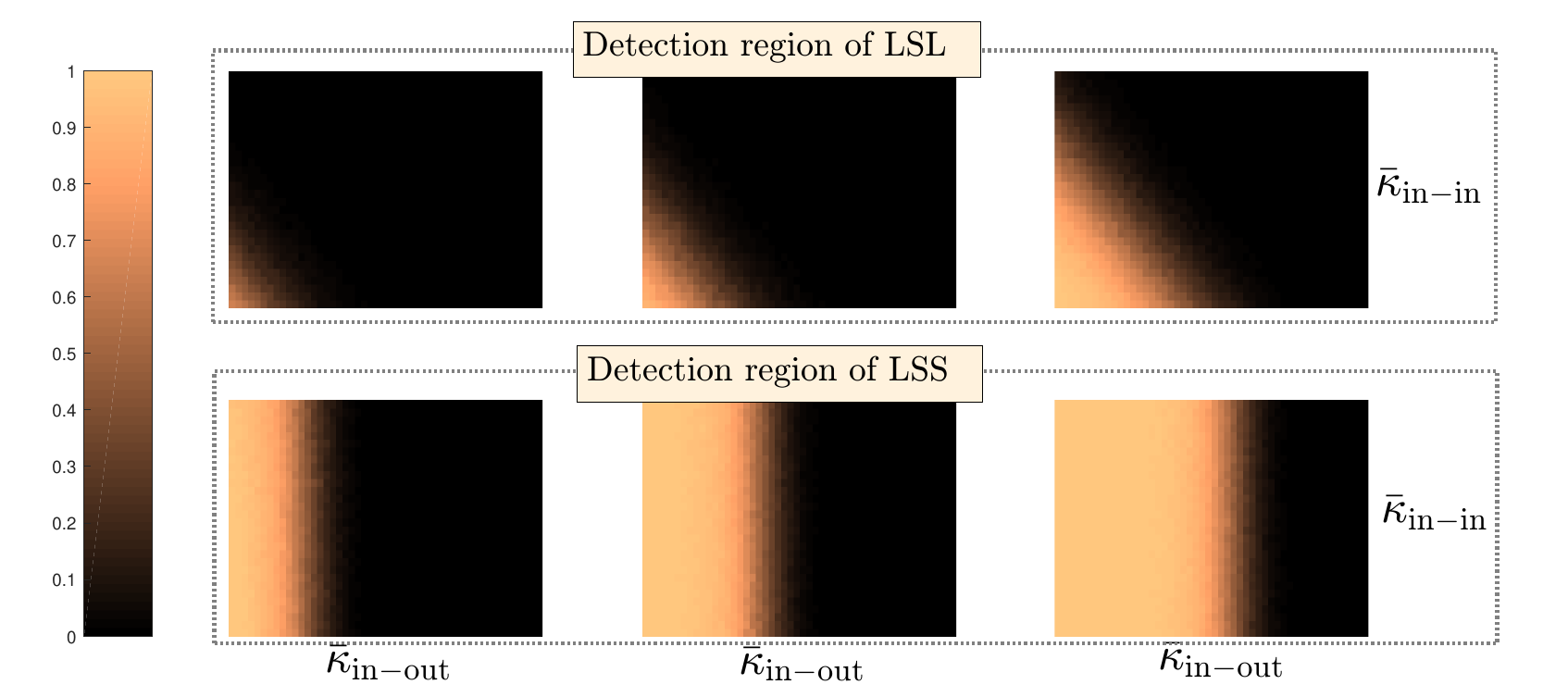}
    \vspace{-13pt}
    \caption{Heatmaps of the error rate
    of the LSL (top row) and the LSS (bottom row)
    estimators in Experiment 2. We chose $n=100$, $m=120$ and 
    $d \in \{10,20,40\}$ from left to right. The 
    parameter $a$ representing the scale of $\kinin$ 
    and corresponding to $y$-axis varies from 0.02 to
    0.08, whereas $b$ representing the scale of 
    $\kinout$ and corresponding to $x$-axis varies 
    from 0.3 to 10. Dark colour means that probability 
    of successful detection is close to 1 (error rate 
    close to zero).}
    \vspace{-10pt}
    \label{fig:heatmap}
\end{figure}

The obtained detection regions are depicted in 
\Cref{fig:heatmap} in the form of heatmaps. This 
visualisation allows us to grasp the forms of the
detection regions for the specific choice of 
parameters considered in this example. The first 
observation is that LSL is clearly superior to 
LSS for all the considered values of the dimension. 
Second, we clearly see the deterioration of the 
detection region when the dimension $d$ becomes
larger. Third, the
values of $\kinout$ used in the plots are at least 
one order of magnitude larger than those of 
$\kinin$. This is in line with the claim of 
\Cref{thm:mild}. We also observe in these pictures
that successful detection occurs when $\kinout$ is 
larger than some threshold even if $\kinin$ is small. 
This must be a nice feature of LSL and LSS in this
specific example, which unfortunately does not
generalize to other examples as shown by our
theoretical results. 

\begin{figure}%
\centering
\includegraphics[width=\textwidth]{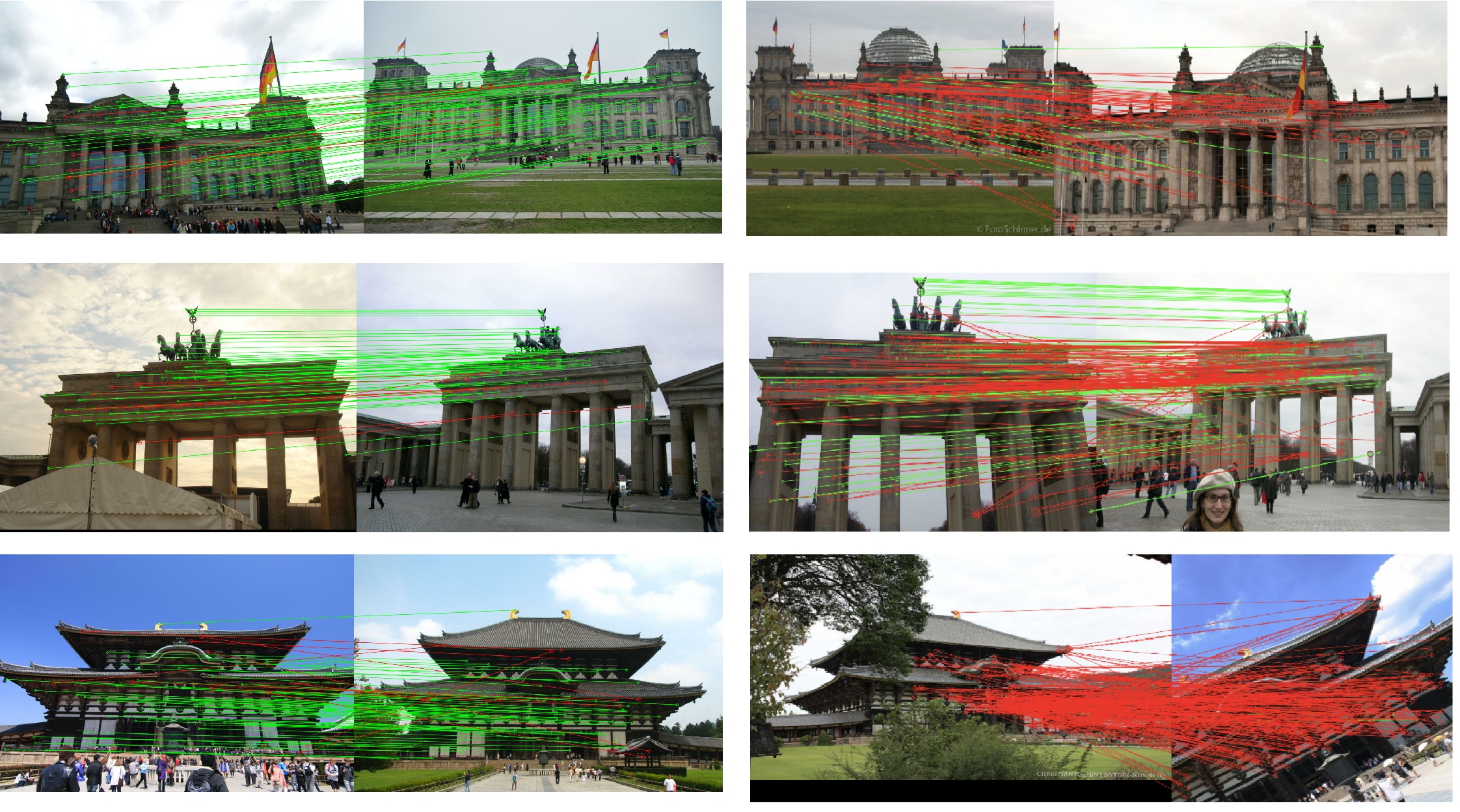}
\caption{Matching map computed by LSL on 
    randomly chosen easy (not challenging) and challenging pairs of images from each scene. 
    The green lines represent the correct 
    matching, and red lines are incorrect ones. 
    }
    \label{fig:samples}
\end{figure}

\paragraph{Experiment 3: 
Real data example}
This experiment is conducted on the IMC-PT 2020 dataset 
from \citep{Jin2020} that consists of images of 16 
different scenes with corresponding 3D point-clouds 
of landmarks, which are used to obtain 
(pseudo) ground-truth local keypoint matchings.  
For a given scene, we sampled 1000 pairs of distinct 
images of the same landmark with different camera 
locations, angles, weather conditions etc. For each 
image pair we generated 2D keypoints from original 
set of 3D points (note that the same 3D point 
appears in both of the images, so we have ground 
truth keypoint matching between 2 images). 
Subsequently, we computed SIFT descriptors 
\citep{lowe2004distinctive} for every keypoint 
in images using Python OpenCV interface \citep{opencv}. 
Some pairs of images being more challenging than others, 
we split the dataset into two sets of image pairs 
in order to gain more understanding on the 
behaviour of the algorithms. The challenging 
pairs are those for which the OpenCV default 
matching algorithm has accuracy less than $50\%$. 
To give a glimpse of what easy and challenging pairs 
of images look like we show in \Cref{fig:samples} 
image pairs with accuracy of OpenCV matching 
algorithm larger than $50\%$ and image pairs with 
accuracy smaller than $50\%$ from each scene.

\begin{figure}%
    \centering
    \includegraphics[width=\textwidth,
    ]{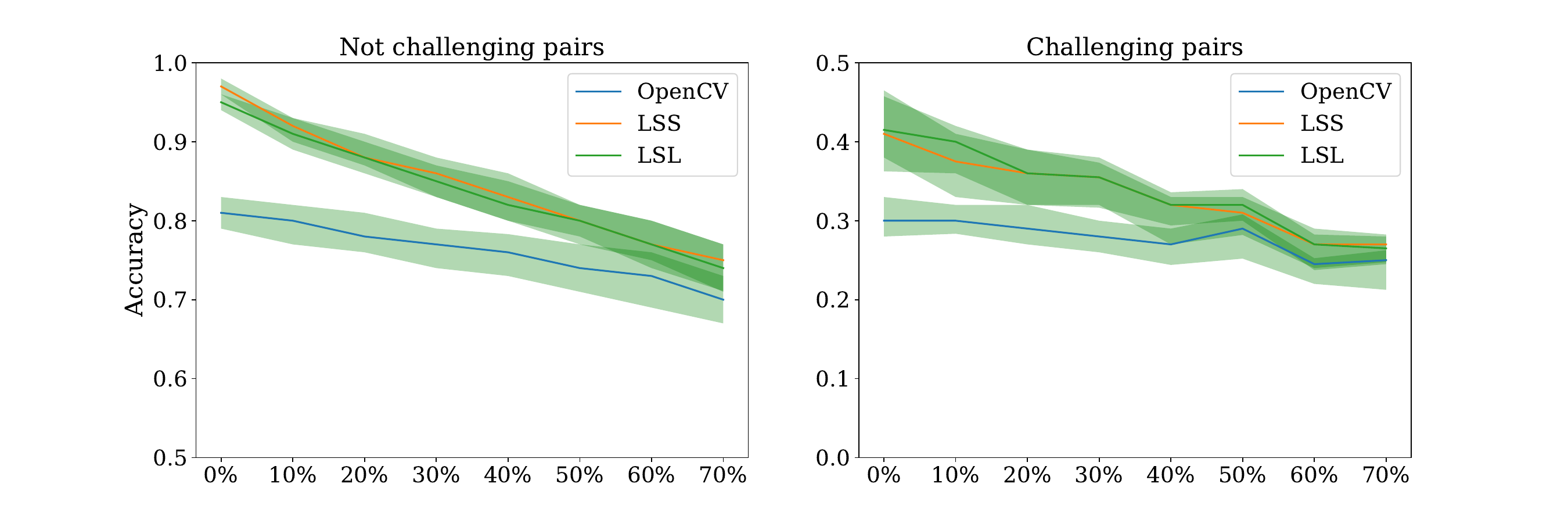}
    \vspace{-10pt}
    
    \caption{The estimation accuracy measured in the Hamming loss of the estimated matching in Exp.~3 for different values of the outlier rate, $(m-n)/n$, varying from $0\%$ to $70\%$. The medians of estimation accuracy both for challenging pairs (right plot) and simple pairs (left plot) of images from Temple Nara  scene was computed using OpenCV, LSS and LSL matchers. The green region represents the interquartile range (lower and upper bounds being $25\%$ and $75\%$ percentiles, respectively).  }
    
    \vspace{-15pt}
    \label{fig:hamming}
\end{figure}

Then, for every image pair, we fixed randomly 
chosen 100 keypoints in the first image (and 
corresponding keypoints in the second image) and 
added outliers to the second image from the 
remaining keypoints. The outlier rate is chosen 
to be between $0\%$ to $70\%$. Finally, 3 
descriptor matching algorithms were applied 
(OpenCV default matching algorithm, LSS and LSL). 
{Note that $\bsigma$ and $\bsigmadiese$ 
from \eqref{model} are unknown, hence LSNS is 
not applicable. One can consider using the 
estimates $\hat\bsigma$ and $\bsigmadiesehat$ 
instead of $\bsigma$ and $\bsigmadiese$ 
in \eqref{LSNS}, but this is beyond the scope 
of this paper.} 

The median estimation accuracy measured in the 
Hamming loss---for the image pairs from Temple 
Nara Japan scene---is plotted in \Cref{fig:hamming}. 
\footnote{We observe very similar behaviour in 
all 3 applied algorithms across other scenes 
as well, therefore the corresponding accuracy 
plots are omitted.} 
The error bars with borderlines corresponding 
to $75\%$ and $25\%$ percentiles are also 
displayed. The first observation is that LSS 
and LSL outperform the OpenCV matcher in terms 
of the number of erroneous matches. Second, the 
rate of correctly matched descriptors deteriorates
with the growth of outlier rate and this deterioration 
seems to be linear. This contrasts with our theoretical
results in which which the impact of the rate of 
outliers is very limited. Note, however, that in the
present experiment the outliers can be very similar 
to the inliers and, therefore, the separation condition 
imposed on $\bar\k_{\textup{in-out}}$ in 
\Cref{thm:upper-bound-LSL,thm:mild} is violated. In 
addition, the results established in this work deal
with the error of detection of matching map and
do not assess the proportion of correctly matched 
descriptors. 

\begin{figure}%
\minipage{0.3\textwidth}
  \includegraphics[width=\linewidth]{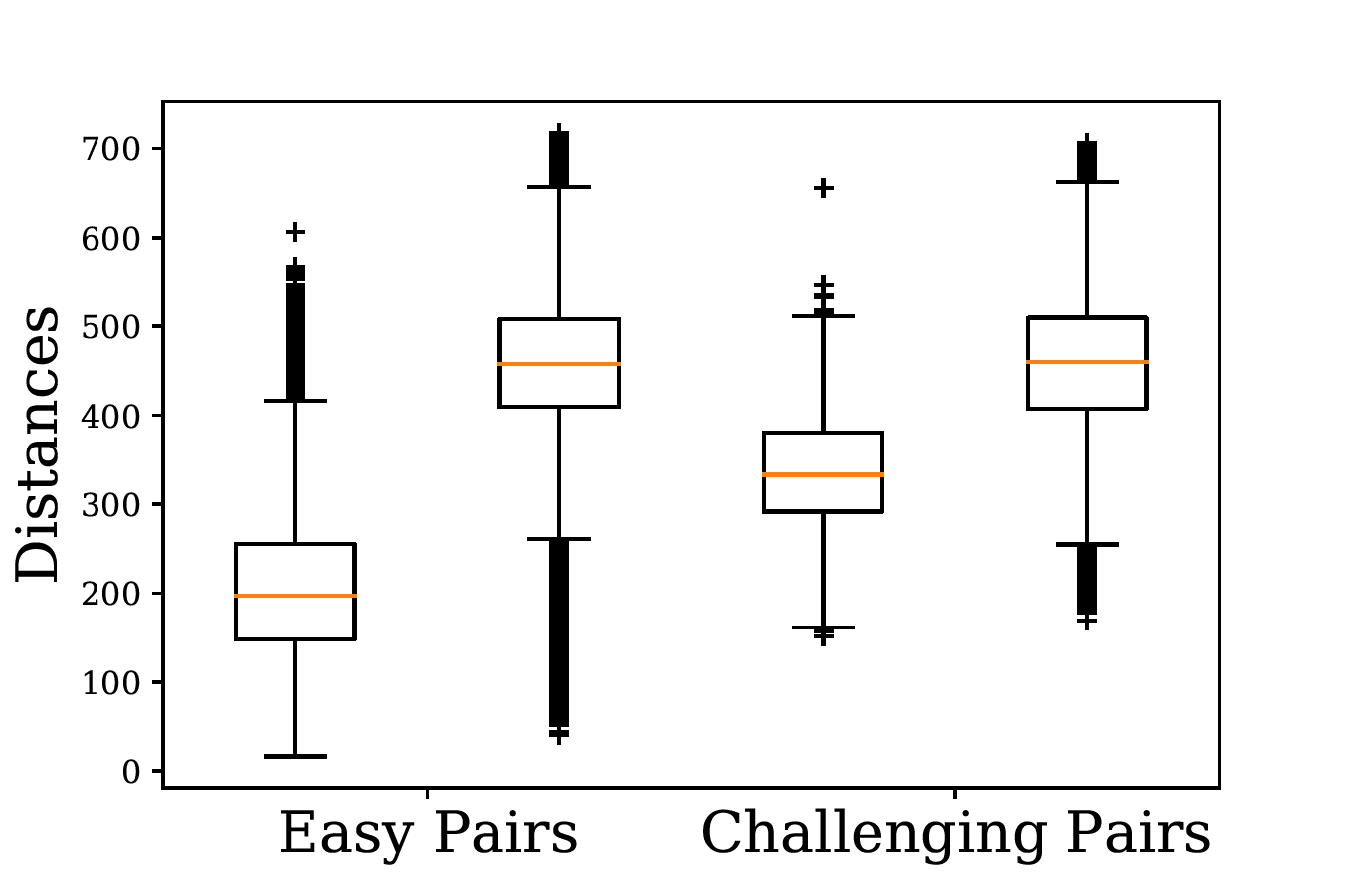}
  
\endminipage\hfill
\minipage{0.3\textwidth}
  \includegraphics[width=\linewidth]{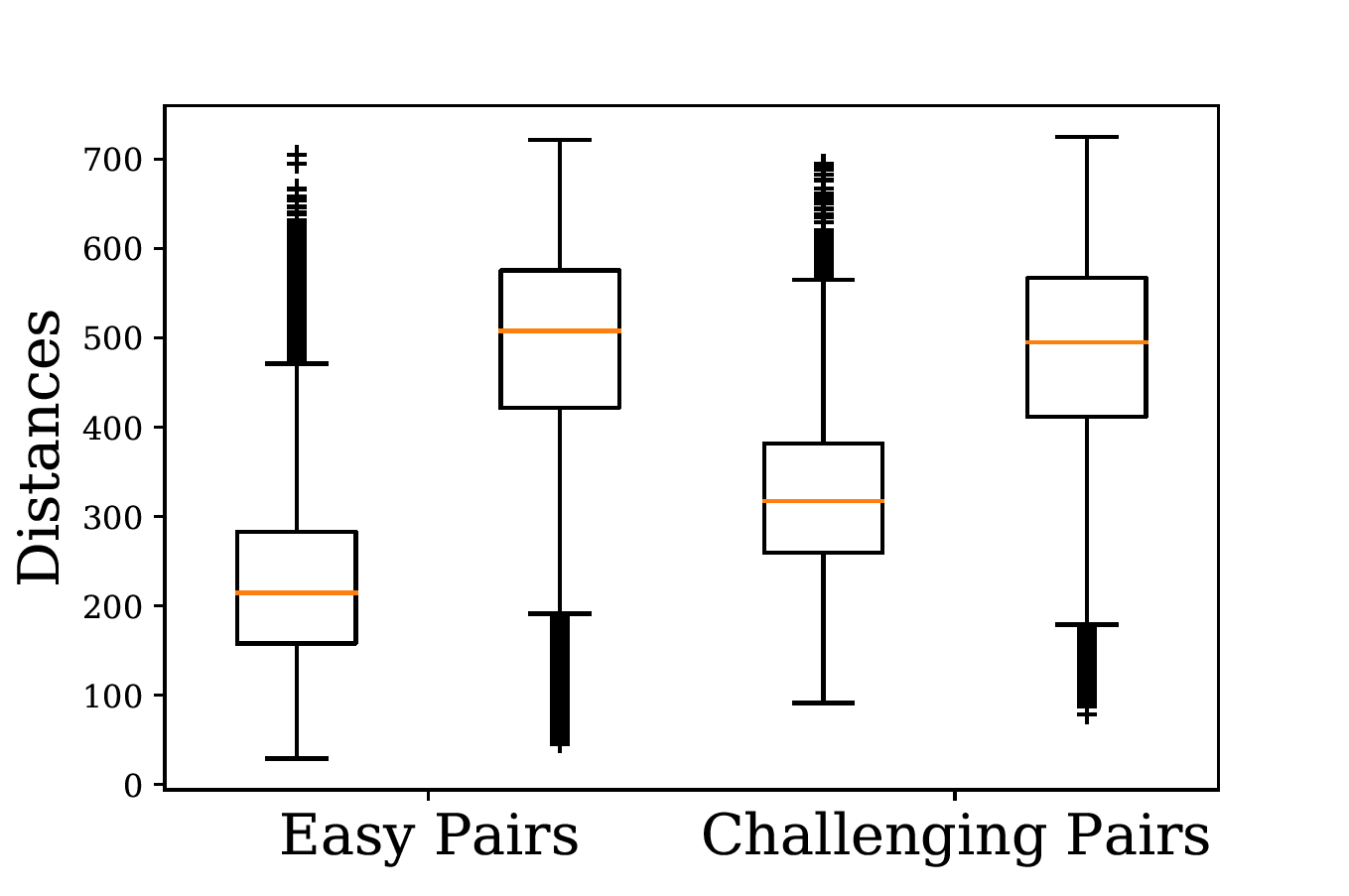}
  
\endminipage\hfill
\minipage{0.3\textwidth}%
  \includegraphics[width=\linewidth]{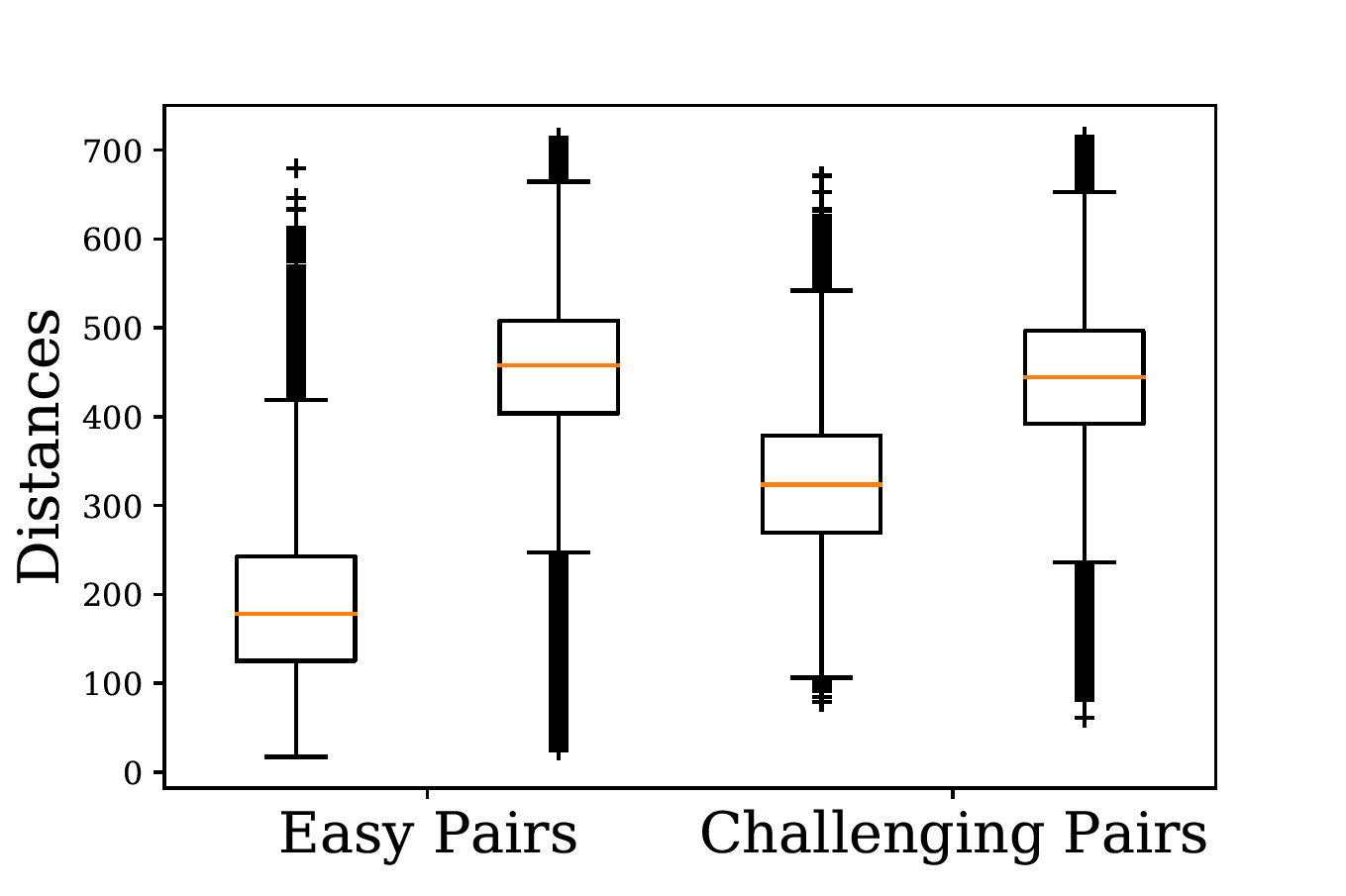}
\endminipage
\caption{ Boxplots of distances between SIFT descriptors 
for  (left to right) Reichstag, Brandenburg Gate and 
Temple Nara scenes. We split datasets into easy and 
challenging pairs according to OpenCV matching algorithm 
score (image pairs with less than $50\%$ correctly matched 
descriptors are considered challenging, the others are easy 
pairs). For each scene we then draw the boxplots of 
distances between descriptors of matching keypoints 
and non-matching keypoints grouped by easy and 
challenging pairs, respectively.  
}
\label{fig:boxplots_dists}
\end{figure}

We also plot the boxplots of the distances between 
SIFT descriptors of matching and non-matching 
keypoints both for easy (not challenging) and 
challenging pairs of images. \Cref{fig:boxplots_dists} 
has 3 plots for each of the scenes and 4 boxplots 
in each of them. The first 2 boxplots correspond 
to the distance between SIFT descriptors of 
matching and non-mathcing keypoints for easy pairs, 
while the last 2 boxplots are that of challenging pairs. 
There are several phenomena that are observed across 
all scenes. First, the median distance for matching 
keypoint descriptors is much smaller than that of 
non-matching keypoint descriptors. Second, the median 
distance between the matching keypoint descriptors 
from challenging pairs is much higher than that of 
easy pairs. We also observe that the distance 
distribution of non-matching keypoint descriptors 
is roughly the same for easy and challenging pairs. 

The results of this experiment suggest 
that LSL and LSS are good estimators in this more
general setting as well (descriptors which are not 
well separated). However, the number of outliers 
might have a significant impact on the accuracy 
of the distance-based algorithms and this impact
needs to be better understood.

\section{Discussion and outlook}\label{app:disc}

\paragraph{Intuitions on the separation rate}
Let us provide some explanations that should help to gain 
intuition on the conditions on $\kinin$ and $\kinout$ obtained 
in our main theorems. More precisely, we will explain in this
paragraph where the right hand side of \eqref{thresh:1} 
comes from. Consider the simpler problem in which we wish 
to test the hypothesis $H_0:\bmu =0$ against $H_1:\bmu\neq 0$ 
based on the observation $\bY$ drawn from the Gaussian 
distribution $\mathcal N_d(\bmu,\sigma^2\mathbf I_d)$. 
This problem has a tight link with the considered problem 
of matching, since one can think of $\bY$ as the difference 
$\bX_i-\bXdiese_j$. We are interested in checking whether
the pair $(i,j)$ is such that $j=\pi^*(i)$, that is whether
$H_0$ is true. 

Using the standard bounds on the tails of the chi-squared 
distribution (\Cref{concentration}), one can check that 
under $H_0$, the random vector $\bY$ lies with probability 
$\ge 1-\alpha$ in the ring $\mathfrak R_0 = B(0,\sigma
\sqrt{d + r_2})\setminus B(0,\sigma\sqrt{d-r_1})$ where 
\begin{align}
    r_1 = 2\sqrt{d\log(1/\alpha)}\qquad\text{and}
    \qquad 
    r_2 = 2\sqrt{d\log(1/\alpha)} + 2\log(1/\alpha). 
\end{align}
Similarly, considering the approximation $\|\bY\|_2^2 
\approx \|\bmu\|_2^2 + \sigma^2\|\bxi\|_2^2$ where $\bxi$ 
is a standard Gaussian vector, we can check that under 
$H_1$, the random vector $\bY$ lies with probability 
$\ge 1-\alpha$ in the ring $\mathfrak R_1 = B(0,\sigma
\sqrt{\|\bmu/\sigma\|_2^2 + d + r_2}) \setminus B(0, 
\sigma\sqrt{\|\bmu/\sigma\|_2^2 + d - r_1})$. 

If the two rings $\mathfrak R_0$ and $\mathfrak 
R_1$ are disjoint, it is possible to decide between 
$H_0$ and $H_1$ by checking whether $\bY$ belongs 
to $\mathfrak R_0$ or not. This condition of disjointness
is equivalent to 
\begin{align}
    \|\bmu/\sigma\|_2^2 + d - r_1 > {d + r_2}.
\end{align}
This leads to 
\begin{align}
    \|\bmu/\sigma\|_2 > \sqrt{r_1 + r_2}  &= \big(4 
    \sqrt{d\log(1/\alpha)} + 2\log(1/\alpha)\big)^{1/2}\\
    &\asymp\big(d\log(1/\alpha)\big)^{1/4}\vee 
    \log^{1/2}(1/\alpha).
\end{align}
The right hand side of the last display is of the same 
order as the right hand side of the \eqref{thresh:1}, for 
small values of $nm$. The fact that for large values of 
$nm$ there is a logarithmic deterioration, due to the 
fact that we have to test a large number of hypotheses 
$H_{0,i,j}: \theta_{\pi^*(i)} = \tdiese_j$, $(i,j)\in 
[n]\times[m]$, is quite common in probability and 
statistics.

\paragraph{Other noise distributions} 
The results of this paper can be extended to sub-Gaussian 
distributions without any change in the rates. The extension 
to sub-exponential distributions seems also possible to do 
using the methodology employed in this paper, but will most 
likely lead to higher-order polylogarithmic terms. 

Finally, considering heavy tailed distributions such as 
the multivariate Student distribution might have stronger 
impact on the rate. Studying this impact is out 
of scope of the present work. 

\paragraph{Outlier detection} The results presented 
in previous sections provide conditions under which 
the objective mapping is identified with high 
probability. This automatically implies that 
the outliers are correctly identified. However, the 
task of outlier detection is arguably simpler than 
that of estimation/detection of $\pi^*$. Therefore, one may 
wonder whether this task can be accomplished under 
weaker assumptions than those required in the 
theorems stated in this paper.  Somewhat surprisingly, 
it turns out that this is not the case unless we 
require the outliers to be very far away from the 
inliers. 

Indeed, on the one hand, if the normalized distance between 
the outliers and the inliers is not larger than $O(d^{1/2})$, 
it follows from the counter-example constructed in the proof
of \Cref{thm:lower-bound} that it is impossible to identify
the outliers using a distance based $M$-estimator. Moreover, this impossibility holds for every estimator $\hat{\pi}$ of the set of outliers, as shown in \Cref{thm:general-lb}.

On the other hand, suitably adapting the arguments of 
the proof of \Cref{thm:mild}, one can prove that if 
the inlier-outlier distance is larger than a threshold 
of order $\sqrt{d}\exp(c n)$ for some $c>0$, the 
LSL recovers the true set of outliers.

\paragraph{Estimation of $\pi^*$ instead of detection}
An interesting yet challenging problem is that of
assessing the minimax risk of estimation of $\pi^*$
when the error is measured, for instance, by means of
the Hamming loss $\ell_{\textup{Hamming}}(\hat\pi;
\pi^*) = \#\{i\in[n] :
\hat\pi(i) \neq \pi^*(i)\}$. It is relevant to study
this problem in a setting where consistent detection
of $\pi^*$ (\textit{i.e.}, Hamming loss equal to zero) 
is impossible, that is when the separation conditions
are violated but some weaker assumptions are satisfied. 
On a related note, one may look for conditions on
the normalized separation distances which ensure the
existence of an estimator $\hat\pi$ such that 
$\mathbf P(\ell_{\textup{Hamming}}(\hat\pi;
\pi^*) \le \tau n)\ge 1-\alpha$. This means that
with probability $\ge 1-\alpha$ the fraction of 
mismatched vectors of the estimated map 
$\hat\pi$ is less than $\tau$, for $\tau\in(0,1)$. 
Note that these problems are not studied even in the
simpler outlier-free framework.

\section{Conclusion}\label{sec:conc}
We have investigated the detection regions in the
problem of detection of the matching map between
two sequences of noisy vectors. We have shown that 
the presence of outliers in one  of the two sequences
has a strong negative impact on the detection 
region. Interestingly, this negative impact is
mitigated in the regime of mild heteroscedasticity, 
\textit{i.e.}, when noise variances are of the same
order of magnitude. In the extremely favorable case 
of homoscedastic noise (all variances are equal), 
the presence of outliers does not make the problem
any harder, provided that the outliers are at least
as different from inliers as two distinct inliers are
different one another. Precise forms of the detection
region in these different cases can be found in
\Cref{tab:1}. The results of the numerical 
experiments conducted on both synthetic and real data
confirm our findings and, furthermore, show the good
behaviour of the LSL estimator in terms of its
robustness to noise and to outliers, not only in 
the problem of detection but also in the problem of 
estimation. In the future, we plan to investigate the
case when both sequences contain outliers and to
obtain theoretical guarantees on the estimation 
error measured by the Hamming distance. 

\appendix
\section{Postponed proofs}\label{app:proofs}

In this appendix we have collected the proofs of the theorems presented in the main text of the paper, as well as some technical definitions used in the proofs. First, denote 
\begin{align}\label{def:kij}
\sigma_{i,j}^2 = \sigma_i^2 + \sdiesesqr_j
\qquad\text{and}\qquad \kappa_{i,j} = 
\frac{\|\theta_i - \tdiese_j\|}{\sigma_{i,j}}
\end{align}
for any pair of indices $ (i,j)$ with $i \in [n]$ 
and $j\in [m]$. We will also use the notation
\begin{align}\label{eq:kappa}
    \bar\k = \min(\kinin,
\kinout).
\end{align}
Second, we define the random variables $\zeta_1$ 
and $\zeta_2$ as follows 
\begin{align}\label{zetas}
\zeta_1 = \max_{i \neq j} \frac{\big|(\theta_i - \tdiese_j)^\top (\sigma_i \xi_i - \sdiese_j \xidiese_j)\big|}{\|\theta_i - \tdiese_j\|\sigma_{i, j}}, \quad \zeta_2 = d^{-1/2} \max_{i, j} \bigg|  \frac{\big\|\s_i \xi_i - \sdiese_j \xidiese_j\big\|^2}{\s_{i,j}^2} - d \bigg|.
\end{align}
It can be easily noticed that $\zeta_1 = \max_{i\neq j} |\zeta_{i,j}|$, where $\zeta_{i,j}$ are standard Gaussian random variables. As for $\zeta_2$, it can be seen that $\zeta_2 = d^{-1/2} \max_{i,j} |\eta_{i,j}|$, where $\eta_{i,j}$ are centered $\chi^2$ random variables with $d$ degrees of freedom, \textit{i.e.} $\eta_{i,j} \stackrel{\mathscr{D}}{=} \chi^2_d - d$. 

In addition, one can infer from \eqref{model}
that for every $i\in[n]$ and every $j\in[m]$, 
we have
\begin{align}
    \| X_i - \Xdiese_j \|^2 &\le \| \t_i - \tdiese_j \|^2 + \sigma_{i, j}^2 (d + \sqrt{d}\, \zeta_2) + 2\zeta_1 \| \t_i - \tdiese_j \| \sigma_{i,j}\\
    & =  \sigma_{i,j}^2\big(\kappa_{i,j}^2 +
    d + \sqrt{d}\, \zeta_2 + 2\zeta_1 \kappa_{i,j} 
    \big),\label{upbX}\\
    \| X_i - \Xdiese_j \|^2 &\ge \| \t_i - \tdiese_j \|^2 + \sigma_{i, j}^2 (d - \sqrt{d}\, \zeta_2) - 2\zeta_1 \| \t_i - \tdiese_j \| \sigma_{i,j}\\
    & =  \sigma_{i,j}^2\big(\kappa_{i,j}^2 +
    d - \sqrt{d}\, \zeta_2 - 2\zeta_1 \kappa_{i,j} 
    \big)\label{lowbX}.
\end{align}

The concentration of the centered and normalized $\chi^2$ random variable, such as $\zeta_2$, is
described in the following lemma.
\begin{lemma}[\citet{LaurentMassart2000}, Eq. (4.3) and (4.4)]\label{concentration}
If $Y$ is drawn from the chi-squared distribution  $\chi^2(D)$, where $D\in\mathbb N^*$, then, for every $x>0$,
\begin{equation}
\begin{cases}
\ \prob\big(Y-D \le -2\sqrt{Dx} \big) \le e^{-x}, \phantom{\Big()}\\
\ \prob\big(Y-D \ge 2\sqrt{Dx} + 2x \big) \le e^{-x}. \phantom{\Big()}
\end{cases}
\end{equation}
As a consequence, for every $y>0$,
$\prob\big(D^{-1/2}|Y-D| \ge y\big) \le 2\exp\big\{-\textstyle\frac18 y(y\wedge \sqrt{D})\big\}$. Or, equivalently, for any 
$\alpha\in(0,1)$, we have
\begin{align}
    \prob\bigg(D^{-1/2}|Y-D| \le 2\sqrt{\log(2/\alpha)} + \frac{2\log(2/\alpha)}{\sqrt{D}}\bigg) \ge 1-\alpha.
\end{align}
\end{lemma}

\subsection{Proof of Theorem~\ref{thm:upperLSNS}}

We prove the upper bound for $\bar\kappa$ in the presence of outliers. Without loss of generality we can assume that  $ \pi^*(i)=i, ~ ~ \forall i \in [n]$. We wish to bound the probability of the event $\Omega = \{ \hat\pi \neq \pi^* \}$, where $\hat\pi = \bar\pi^{\text{LSNS}}$. It is evident that 
\begin{align}\label{omega1}
\Omega \subset \bigcup_{\pi\neq\pi^*} \Omega_{\pi},
\end{align}
where the union is taken over all possible 
injective mappings $\pi:[n]\to[m]$ and 
\begin{align}
\Omega_{\pi} &= \bigg\{ \sum_{i=1}^n \frac{\| X_i - \Xdiese_i \|^2}{2\sigma_i^2} \ge \sum_{i=1}^n \frac{\| X_i - \Xdiese_{\pi(i)}\|^2}{\sigma_i^2 + (\sdiese_{\pi(i)})^2 }\bigg\}.
\end{align}
One easily checks that the following inclusion holds:
\begin{align}\label{omega2}
\Omega_{\pi}  \subset
\bigcup_{i=1}^n \bigcup_{j \in [m] \setminus \{ i \}} \bigg\{ \frac{\| X_i - \Xdiese_i \|^2}{2\sigma_i^2} \ge \frac{\| X_i - \Xdiese_{j}\|^2}{\sigma_i^2 + (\sdiese_{j})^2 }\bigg\}.
\end{align}
Since $\pi^*(i) = i$ for every $i \in [n]$, 
$\kappa_{i,i} = 0$ (see the definition in 
\eqref{def:kij}) and, in view of \eqref{upbX}, 
\begin{align}\label{eq:i=j}
\| X_i - \Xdiese_i \|^2  
\le 2\s_i^2 (d + \sqrt{d}\, \zeta_2).
\end{align}
Similarly, for every $j \in [m]$ and $j \neq i$, 
in view of \eqref{lowbX},
\begin{align}
\| X_i - \Xdiese_j \|^2 &\ge 
\s_{i,j}^2 (\kappa_{i,j}^2 + d - \sqrt{d}\,\zeta_2 
- 2 \kappa_{i,j} \zeta_1).
\end{align}
Recall that $\bar\kappa$ defined in \eqref{eq:kappa},
is the smallest normalized distance $\kappa_{i,j}$.
Therefore, on the event $\Omega_1 = \{ \bar\k \ge
\zeta_1\}$, the previous display implies that 
\begin{align}\label{eq:ineqj}
\frac{\| X_i - \Xdiese_j \|^2}{\s_{i,j}^2} \ge \bar\k^2 - 2\bar\k \zeta_1 + d - \sqrt{d}\, \zeta_2.
\end{align}
Hence, combining obtained bounds \eqref{eq:i=j} and \eqref{eq:ineqj} we get that
\begin{align}
    \bigg\{ \frac{\| X_i - \Xdiese_i \|^2}{2\sigma_i^2} \ge \frac{\| X_i - \Xdiese_{j}\|^2}{\sigma_i^2 + (\sdiese_{j})^2 }\bigg\}\cap \Omega_1 &\subset
    \Big\{ d + \sqrt{d}\, \zeta_2 \ge
    \bar\k^2 - 
    2\bar\k\zeta_1 + d - 
    \sqrt{d}\, \zeta_2\Big\}\\
    & = \Big\{ 2\sqrt{d}\, \zeta_2  + 
    2\bar\k\,\zeta_1 \ge
    \bar\k^2  \Big\}.\label{omega3}
\end{align}
Note that the event on the right hand side of the
last display is independent of the pair $(i,j)$. 
This implies that
\begin{align}
\Omega \cap \Omega_1 
    &\stackrel{\textup{by \eqref{omega1}}}{\subset} 
    \bigg(\bigcup_{\pi\neq\pi^*} \Omega_{\pi} 
    \bigg)\cap \Omega_1\\
    &\stackrel{\textup{by \eqref{omega2}}}{\subset}   \bigg(\bigcup_{i=1}^n \bigcup_{j \in [m] \setminus \{ i \}} \bigg\{ \frac{\| X_i - \Xdiese_i \|^2}{2\sigma_i^2} \ge \frac{\| X_i - \Xdiese_{j}\|^2}{\sigma_i^2 + (\sdiese_{j})^2 }\bigg\} \bigg)\cap \Omega_1\\
    & \stackrel{\hphantom{\textup{by \eqref{omega3}}}}{\subset} \bigcup_{i=1}^n \bigcup_{j \in [m] \setminus \{ i \}} \bigg(\bigg\{ \frac{\| X_i - \Xdiese_i \|^2}{2\sigma_i^2} \ge \frac{\| X_i - \Xdiese_{j}\|^2}{\sigma_i^2 + (\sdiese_{j})^2 }\bigg\} \cap \Omega_1 \bigg)\\
    & \stackrel{\textup{by \eqref{omega3}}}{\subset}  \Big\{ 2\sqrt{d}\, \zeta_2  + 
    2\bar\k\,\zeta_1 \ge
    \bar\k^2  \Big\}. \label{eq:omega_omega1}    
\end{align}
Using \eqref{eq:omega_omega1} we can show that
\begin{align}
\prob(\O)
        &\le \prob(\O_1^\complement) + \prob\big(\O\cap\O_1\big)\\
        &\le \prob\big(\z_1\ge \bar\k\big)
            +\prob(2\sqrt{d}\z_2+2\bar\k\z_1\ge\bar\k^2)\\
        &\le \prob\big(\z_1\ge \bar\k\big) + 
            \prob\big(\z_1\ge {\textstyle\frac14}\bar\k\big)
            +\prob\Big(2\sqrt{d}\z_2+2\bar\k\z_1\ge\bar\k^2; \z_1< {\textstyle\frac14}\bar\k\Big)\\
        &\le 2\prob\big(\z_1\ge {\textstyle\frac14}\, {\bar\k}\big)+
            \prob\bigg(\z_2\ge \frac{\bar\k^2}{4\sqrt{d}}\bigg).
            \label{ineq:2}
\end{align}
For suitably chosen standard Gaussian random variables $\z_{i,j}$ 
it holds that $\z_1 = \max_{i\not=j}|\z_{i,j}|$. Therefore, using the tail bound for the standard Gaussian
distribution and the union bound, we get
\begin{align}
\prob\Big(\z_1\ge {\textstyle\frac14}\bar\k\Big)
    &\le \sum_{i\not =j} \prob\Big(|\z_{i,j}|\ge \textstyle\frac14\bar\k\Big)
    \le 2 n m\, e^{-\bar\k^2/32}.\label{ineq:3}
\end{align}
To complete the proof, it remains to upper bound
the second term in the right hand side of \eqref{ineq:2}, \textit{i.e.}, to evaluate the
tail of the random variable $\zeta_2$. To this 
end, we use the concentration result stated in \Cref{concentration} with $y=\frac{\bar\k^2}{4\sqrt{d}}$, combined with the union bound and simple algebra. This yields 
\begin{align}
\prob\Big(\z_2 \ge \frac{\bar\k^2}{4\sqrt{d}}\Big)
        &\le 2 n m \exp\Big\{-\frac 1 8 \cdot  \frac{\bar\k^2}{4\sqrt{d}}\Big(\frac{\bar\k^2}{4\sqrt{d}} \wedge \sqrt{d}\Big)\Big\} \\
        &= 2 n m \exp\Big\{-\frac{(\bar\k/16)^2}{d}(2\bar\k^2\wedge 8d)\Big\},\label{ineq:4}
\end{align}
where the $nm$ factor in front of the exponent comes from the union bound for all $nm$ pairs $(i, j)$ from the definition of $\zeta_2$, while the exponent is a direct application of \Cref{concentration}. Finally, using inequalities (\ref{ineq:2})-(\ref{ineq:4}), we get that whenever 
\begin{align}\label{kappa-lsns}
\bar\k\ge 4\Big(\sqrt{2\log(8n m/\a)} \vee \big(d\log(4n m/\a)\big)^{1/4}\Big),
\end{align}
the probability of incorrect matching is at most $\alpha$. Thus, we have formally showed that if \eqref{kappa-lsns} holds then $\prob(\hat\pi\neq \pi^*)=\prob(\O)\le \a,$ as desired.

\subsection{Proof of Theorem \ref{thm:upper-bound-LSL}}
We prove the upper bound for $\bar\kappa = \min(\kinin,
\kinout)$ in the presence of outliers and in the case of unknown noise variance. We wish to bound the probability of the event $\Omega = \{\hat\pi \neq \pi^*\}$, where $\hat\pi = \hat\pi^{\textup{LSL}}$ and $\pi^*(i) = i$ for all $i\in[n]$. It is evident that 
\begin{align}\label{omega-in-lsl}
    \Omega \in \bigcup_{\pi \neq \pi^*} \Omega_{\pi},
\end{align}
where 
\begin{align}
    \Omega_{\pi} &= \Big\{\sum_{i=1}^n \log \| X_i - \Xdiese_i\|^2 \ge \sum_{i=1}^n \log \| X_i - \Xdiese_{\pi(i)}\|^2\Big\} \\ 
    &\subset \bigcup_{i=1}^n \bigcup_{j \in [m] \setminus \{i\}} \Big\{ \log \| X_i - \Xdiese_i \|^2 \ge \log \| X_i -\Xdiese_j\|^2\Big\}
    \label{omega2-in-lsl}
\end{align}
Recall that $\bar\k = \min(\kinin,
\kinout)$. On the event $\Omega_1 = \{ \bar\kappa \ge \zeta_1 \}$, from equation \eqref{lowbX}, 
we get  
\begin{align}\label{lwb-X}
    \frac{\| X_i - \Xdiese_j \|^2}{\sigma_{i, j}^2} \ge \bar\kappa^2 - 2\zeta_1 \bar\kappa + d - \sqrt{d}\, \zeta_2.
\end{align}
Note that the expression on the right hand side 
of the last display is independent of the 
pair $(i,j)$.  This implies that
\begin{align}
\Omega \cap \Omega_1 
    &{\subset} 
    \bigg(\bigcup_{\pi\neq\pi^*} \Omega_{\pi} 
    \bigg)\cap \Omega_1  \tag*{[\textup{by \eqref{omega-in-lsl}}]}\\
    &{\subset}   \bigg(\bigcup_{i=1}^n \bigcup_{j \in [m] \setminus \{ i \}} \bigg\{ \log \| X_i - \Xdiese_i \|^2 \ge \log \| X_i -\Xdiese_j\|^2\bigg\} \bigg)\cap \Omega_1 \tag*{[\textup{by \eqref{omega2-in-lsl}}]}\\
    & {\subset} \bigcup_{i=1}^n \bigcup_{j \in [m] \setminus \{ i \}} \bigg(\bigg\{  \| X_i - \Xdiese_i \|^2 \ge  \| X_i -\Xdiese_j\|^2\bigg\} \cap \Omega_1 \bigg)  \\
    & {\subset}  \Big\{ 2\sigma_i^2 (d+\sqrt{d}\zeta_2) \ge \sigma_{i, j}^2 (\bar\kappa^2 - 2\zeta_1 \bar\kappa + d - \sqrt{d}\, \zeta_2) \Big\}  \tag*{[\textup{by \eqref{upbX},\eqref{lwb-X}}]}\\
    & {\subset}  \Big\{ 2 (d+\sqrt{d}\zeta_2) \ge \bar\kappa^2 - 2\zeta_1 \bar\kappa + d - \sqrt{d}\, \zeta_2 \Big\},  \tag*{[\textup{since } $\sigma_i\le \sigma_{i,j}$]}\\
    & {\subset}  \Big\{ 3\sqrt{d}\,\zeta_2 + 2 \zeta_1\bar\kappa \ge \bar\kappa^2 - d \Big\}. \label{eq:omega-event-lsl}    
\end{align}
We can bound the probability
of incorrect matching $\prob(\O)$ using the relationship obtained in \eqref{eq:omega-event-lsl}
\begin{align}
    \prob(\Omega) &\le \prob(\Omega_1^{\complement}) + \prob(\Omega \cap \Omega_1) \\
    &\le \prob(\zeta_1 \ge \bar\kappa) + \prob\big(3\sqrt{d}\,\zeta_2 + 2 \zeta_1 \bar\kappa \ge  \bar\kappa^2 - d\big).
\end{align}
From the last inequality, we infer that
\begin{align}
    \prob(\Omega) 
    &\le \prob(\zeta_1 \ge \bar\kappa) + \prob\Big(\zeta_1 \ge {\textstyle\frac14} {\bar\kappa}\Big) + \prob\Big(3\sqrt{d}\,\zeta_2 + 2 \zeta_1 \bar\kappa \ge  \bar\kappa^2 - d\,;\,
    \zeta_1 < {\textstyle\frac14} \bar\kappa\Big) \\
    &\le 2\prob\Big(\zeta_1 \ge {\textstyle\frac14} {\bar\kappa}\Big) + \prob\Big(3\sqrt{d}\,\zeta_2 \ge {\textstyle\frac12\bar\kappa^2} - d\Big)\\
    &\le 2\prob\Big(\zeta_1 \ge {\textstyle\frac14} {\bar\kappa}\Big) + \prob\Big(\zeta_2 \ge \frac{\bar\kappa^2 - 2d}{6\sqrt{d}}\Big).
    \label{prob-error}
\end{align}
As mentioned in the beginning of the section, for 
suitably chosen standard Gaussian random variables 
$\z_{i,j}$ it holds that $\z_1 =
\max_{i\not=j}|\z_{i,j}|$. Therefore, using the tail 
bound for the standard Gaussian
distribution and the union bound, we get
\begin{align}\label{z1-bound}
\prob\Big(\z_1\ge {\textstyle\frac14}\bar\k\Big)
    &\le \sum_{i\not =j} \prob\Big(|\z_{i,j}|\ge \textstyle\frac14\bar\k\Big)
    \le 2 n m\, e^{-\bar\k^2/32}\le \alpha/4.
\end{align}
To complete the proof, it remains to upper bound
the second term in the right hand side of \eqref{prob-error}, \textit{i.e.}, to evaluate the
tail of the random variable $\zeta_2$. Using \Cref{concentration} with $y =  (\bar\kappa^2 -2d)/(6\sqrt{d})$---which is positive under the conditions of the theorem---combined with the union bound, we arrive at 
\begin{align}
    \prob(\zeta_2 \ge y) &\le 2nm \exp\Big\{- \frac 1 8 y(y\wedge \sqrt{d})\Big\}\\
    & = 2nm \bigg(\exp\Big\{- \frac 1 8 y^2\Big\}\vee \exp\Big\{- \frac 1 8 y\sqrt{d}\Big\}\bigg).
\end{align}
One easily checks that the last expression is smaller than 
$\alpha/2$ if and only if
\begin{align}
    y^2\ge 8 \log(4nm/\alpha)\qquad\text{and}\qquad
    y\sqrt{d}\ge 8 \log(4nm/\alpha)
\end{align}
which is equivalent to
\begin{align}
    y\ge \big(2 \sqrt{2\log(4nm/\alpha)}\big)\vee 
    \big( (8/\sqrt{d}) \log(4nm/\alpha)\big).
\end{align}
Replacing $y =  (\bar\kappa^2 -2d)/(6\sqrt{d})$, the last
inequality becomes
\begin{align}
    \bar\kappa^2\ge 2d + \big(12 \sqrt{2d\log(4nm/\alpha)}\big)\vee 
    \big( 48 \log(4nm/\alpha)\big).
\end{align}
Combining the inequality from the last display with the bound derived from \eqref{z1-bound} we get that all these bounds are satisfied whenever 
\begin{align}
    \bar\kappa \ge \sqrt{2d} + 4\bigg\{\bigg(2d\log \frac{4nm}{\alpha}\bigg)^{1/4} \vee \bigg(3\log\frac{8nm}{\alpha}\bigg)^{1/2}\bigg\}.
\end{align}
Therefore, under this condition on $\bar\kappa$, the probability of the incorrect matching is at most $\alpha$, \textit{i.e.} $\prob(\hat\pi \neq \pi^*) = \prob(\O) \le \alpha$. 

\subsection{Proof of Theorem \ref{thm:lower-bound}}
First we fix $m = n + 1 $ and $\pi^*(i) = i$ for all $i \in [n]$, where  $\pi^*$ is the correct matching. Let $\sdiese_1 = 1$ and $\sdiese_{i+1} = \alpha^i$ for all $i \in [n]$, where $\alpha \ll 1$. Then let's take $\pi(i) = i + 1$ for all $i \in [n]$. Let $L(\pi)$ be the vector of distances $\| X_i - \Xdiese_{\pi(i)}\|$ for a matching scheme $\pi$
\begin{align}
    L(\pi) = 
    \begin{bmatrix} 
        \| X_1 - \Xdiese_{\pi(1)}\| \\ 
        \| X_2 - \Xdiese_{\pi(2)}\| \\ 
        \cdots 
        \\ 
        \| X_n - \Xdiese_{\pi(n)}\|
    \end{bmatrix}.
\end{align}
The next lemma shows that the event $L(\bar\pi) < L(\pi^*)$ (coordinate-wise) occurs with probability at least $1/4$. 

\begin{lemma}\label{lem:counterexample}
Let $n \ge 4$, $d \ge 422\log(4n)$ and 
$\tdiese_1 = (1;0;\ldots;0)^\top$. Assume that
$\pi^*(i)=i$, $\sdiese_i = 2^{-(i-1)}$ and 
$\tdiese_{i+1} = \tdiese_i + 2^{-(i+1)}\sqrt{d} \,
\tdiese_1$ for all $i \in [n+1]$. Then $L(\pi^*) > L(\bar\pi)$ with probability greater than $1/4$, where $\bar\pi$ is the injection
defined by $\bar\pi(i) = i+1$. Furthermore, for
these values $(\btdiese,\bsigmadiese,\pi^*)$, we 
have $\k_{\textup{in-in}} = \k_{\textup{in-out}} 
= \sqrt{d/20}$.
\end{lemma}
\begin{proofof}[Proof of  \Cref{lem:counterexample}]
Let us denote 
\begin{align}
    \bar\kappa_i \triangleq \frac{\| \tdiese_{\pi(i)} - \theta_i\|}{\sqrt{\sigma_i^2 + \sdiesesqr_{\bar\pi(i)}}} = \sqrt{d/20}, 
    \qquad \text{for all } i \in [n].
\end{align}
Recall that $\sigma_{i,j}^2 = \sigma_i^2 + \sdiesesqr_j$ and write 
\begin{align}
    L_i(\pi) = \| X_i - \Xdiese_{\pi(i)} \|^2 &= \| \t_i - \tdiese_{\pi(i)} + \zeta_i \sigma_{i, \pi(i)} \|^2,
\end{align}
where $\zeta_i \sim \mathcal{N}(0, I_d)$. Notice that $L_i(\pi^*) = 2\sigma_i^2 \| \zeta_i \|^2$ for all $i \in [n]$. Similarly, the expression from the last display for $\bar\pi$ reads as 
\begin{align}
    L_i(\bar\pi) = \| \zeta_i \sigma_{i, \bar\pi(i)}\|^2 \Big( 1 + \frac{\bar\kappa_i^2}{\| \zeta_i \|^2}\Big) + {2\sigma_{i, \bar\pi(i)}\zeta_i^\top (\theta_i - \tdiese_{\bar\pi(i)})}.
\end{align}
Plugging in the values of $\bsigmadiese$ with $\alpha = 1/2$ and $\bar\pi(i) = i+1$ we arrive at 
\begin{align}
    L_i(\pi^*) = 2^{3-2i} \| \zeta_i\|^2, \qquad L_i(\bar\pi) = \frac{5}{2^{2i}} \| \bar\zeta_i \|^2\Big(1  + \frac{\bar\kappa_i^2}{\| \bar\zeta_i\|^2}\Big) + \frac{\sqrt{5}}{2^{i-1}} \bar\zeta_i^\top(\theta_i - \tdiese_{i+1}),
\end{align}
where in the second expression we write $\bar\zeta_i$ instead of $\zeta_i$ to indicate that these random variables are different, though both are standard normal $d$-dimensional vectors. We first replace the second term of $L_i(\bar\pi)$ with its upper bound that holds with probability of at least $1/4$. It is evident that the random variable $Z \triangleq {2\sigma_{i, \bar\pi(i)}\zeta_i^\top (\theta_i - \tdiese_{\bar\pi(i)})}$ is Gaussian with standard deviation $\sigma \triangleq 2\sigma_{i, \bar\pi(i)} \| \theta_i - \tdiese_{\bar\pi(i)}\| = 2\sigma_{i, \bar\pi(i)}^2 \bar\kappa_i$, therefore
\begin{align}
    \prob(Z \ge \sigma \sqrt{2\log 4}) \le \frac 1 4.
\end{align}
Hence, on the event $\Omega = \{ Z \le 2\sigma_{i, \bar\pi(i)}^2 \bar\kappa_i \sqrt{2\log 4}\}$ the inequality $L_i(\pi^*) > L_i(\bar\pi)$ holds whenever 
\begin{align}
    &\frac 8 {2^{2i}} \| \zeta_i\|^2 > \frac{5}{2^{2i}} \| \bar\zeta_i\|^2\Big(1  + \frac{\bar\kappa_i^2}{\| \bar\zeta_i\|^2}\Big) + \frac{5}{2^{2i}} \bar\kappa_i \sqrt{8\log 4},\nonumber \\
    &\frac 8 5 \| \zeta_i \|^2 - \| \bar\zeta_i \|^2 > \bar\kappa_i^2 + 2\bar\kappa_i\sqrt{2\log 4} \label{chi2-diff}.
\end{align}
Notice that the left hand side of \eqref{chi2-diff} is a weighted difference of two centered and normalized $\chi^2$ random variables with $d$ degrees of freedom. The concentration inequality for such difference is a direct consequence of \Cref{concentration}. Namely, for $X, Y \sim \chi^2_d$ the concentration bound for $Z = \alpha X - \beta Y$ with arbitrary $\alpha, \beta \in \mathbb{R}$ reads as  
\begin{align}
    \prob(Z \ge (\alpha - \beta)d - 2\sqrt{dx}(\alpha + \beta) - 2\beta x) \ge 1 - 2 e^{-x}.
\end{align}
It is easy to verify that given $n \ge 4, d \ge 422 \log (4n)$ and $\bar\kappa_i \le \sqrt{d/20}$, then 
\begin{align}
    \bar\kappa_i^2 + 2 \bar\kappa_i \sqrt{2 \log 4} \le \frac 3 5 d - \frac{26}{5} \sqrt{d\log(4n)} - 2\log(4n),
\end{align}
where the right hand side is the quantile of $Z$ with $x = \log(4n)$.
Combining the inequality from the last display with \eqref{chi2-diff} we get that on the event $\Omega$ we have
\begin{align}
    \prob(L_i(\pi^*) > L_i(\bar\pi)) \ge 1 - \frac{1}{2n}. 
\end{align}
Recall that $\prob(\Omega) \ge 3/4$, then using the union bound for events $\Omega$ and $\{L_i(\pi^*) > L_i(\bar\pi)\}$ all $i \in [n]$ we arrive at  $\prob(L(\pi^*) > L(\bar\pi)) > 1/4$. This completes
the proof of \Cref{lem:counterexample}.

\end{proofof}

Therefore, using the result of \Cref{lem:counterexample} and applying any non-decreasing function $\rho(\cdot)$ to each of the coordinates of $L(\bar\pi)$ and $L(\pi^*)$ yields 
\begin{align}
    \sum_{i=1}^n \rho_i(\|X_i - \Xdiese_{\bar\pi(i)}\|) < \sum_{i=1}^n \rho_i(\|X_i - \Xdiese_{\pi^*(i)}\|)
\end{align}
with probability of at least $1/4$. This, in turn, implies that an optimizer will not choose $\pi^*$ 
on this event. Hence,  $\prob(\bar\pi \neq \pi^*) > \nicefrac 1 4$, concluding the proof of the theorem. 

\subsection{Proof of \Cref{thm:general-lb}}
We denote the set of all injective functions $\pi :[n] \to [m]$ as $\mathfrak{I}_{n,m}$. We use the notation $D(\prob, \bQ)$ for the Kullback-Leibler (KL) divergence between two probability measures $\prob$ and $\bQ$ such that $\prob$ is absolutely continuous with respect to $\bQ$, $\prob \ll \bQ$. The identity mapping denoted by \textit{id} is defined as follows: $id(i) = i, \, \, \forall i \in [n]$. It is also assumed that $\pi^* = id$.

To establish the general lower bound we use the following lemma: 
\begin{lemma}[\citet{Tsybakov2009}, Theorem 2.5]\label{lower-bound}
    Assume that for some integer $M\ge 2$ there exist distinct injective functions $\pi_0,\ldots,\pi_M \in \mathfrak{I}_{n, m}$ and mutually absolutely continuous probability measures $\bQ_0,\ldots,\bQ_M$ defined on a common probability space $(\mathcal Z,\mathscr Z)$ such that
    \begin{align}
        \frac 1 M \sum_{j=1}^M D(\bQ_j, \bQ_0) \le \frac 1 8 \log M.
    \end{align}
    Then, for every measurable mapping $\tilde{\pi} : \mathcal{Z} \to \mathfrak{I}_{n, m}$, 
    \begin{align}
        \max_{j=0, \dots, M} \bQ_j(\tilde{\pi} \neq \pi_j) \ge \frac{\sqrt{M}}{\sqrt{M} + 1}\Big(\frac 3 4 - \frac{1}{2\sqrt{\log(M)}}\Big).
    \end{align}
\end{lemma}

Since $d \ge 16 \log(nm)$ then the rate from \Cref{thm:upperLSNS} becomes of order $(d\log(nm))^{1/4}$. We show that for $6\kappa \ge (d \log(nm))^{1/4}$ there is indeed a setting where the detection of $\pi^*$ fails with probability at least $\nicefrac{1}{4}$ for any matching map $\tilde{\pi}\in\mathfrak{I}_{n,m}$. To show this we use \Cref{lower-bound} with properly chosen family of probability measures described in the following lemma. 

\begin{lemma}[\citet{collier2016minimax}, Lemma 14]
Let $\varepsilon_1, \dots, \varepsilon_m$ be real numbers defined by 
\begin{align}
    \varepsilon_k = \sqrt{2/d} \, \kappa \sdiese_k, \quad \forall k \in [m],
\end{align}
and let $\mu$ be the uniform distribution on $\mathcal{E} = \{ \pm \varepsilon_1\}^d \times \dots \times \{ \pm \varepsilon_m \}^d$. Denote by $\prob_{\mu, \pi}$ the probability measure on $\mathbb{R}^{d\times m}$ defined by $\prob_{\mu, \pi}(A) = \int_{\mathcal{E}} \prob_{\btheta, \pi}(A) \mu(d\btheta)$. Let $\bar{\Theta}_{\kappa}$ be the set of $\btdiese$ such that $6\kappa \ge (d\log(nm))^{1/4}$. Assume that $\sdiese_1 \le \dots \le \sdiese_m$ and $\sdiesesqr_m / \sdiesesqr_1 \le 1 + \sqrt{\frac{\log(nm)}{16d}}$. Let $\pi = (k \, \, k')$ be the transposition that only permutes $k^{\textup{th}}$ and $k'^{\textup{th}}$ observations ($k < k')$. Then, the Kullback-Leibler divergence between $\prob_{\mu, \pi}$ and $ \prob_{\mu, id}$ can be bounded as follows
\begin{align}
    D(\prob_{\mu, \pi}, \prob_{\mu, id}) \le \frac 1 8 \log(m(m-1)/2).
\end{align}
Additionally, $\mu(\mathcal{E} \setminus \bar\Theta_{\kappa}) \le (m(m-1)/2) e^{-d/8}$.
\end{lemma}

Applying \Cref{lower-bound} with $M = m(m-1)/ 2$, $\bQ_0 = \prob_{\mu, id}$ and $\{\bQ_j\}_{j=1, \dots, M} = \{ \prob_{\mu, \pi_{k, k'}}\}_{k\neq k'}$ we obtain that for any estimator 
$\hat\pi$
\begin{align}
 \max_{\pi^* \in \mathfrak{I}_{n, m}} \sup_{\btdiese \in \bar{\Theta}_{\kappa}} \prob_{\btdiese, \bsigmadiese, \pi^*}(\hat\pi \neq \pi^*) &\ge \max_{\pi^*\in\{\id\}\cup\{\pi_{k,k'}\}}
\int_{\bar\T_\k} \prob_{\bthetadiese,\pi^*} \big( \hat\pi \neq \pi^* \big) \frac{\mu(d\bthetadiese)}{\mu(\bar\T_\k)} \\
&\ge \max_{\pi^*\in\{\id\}\cup\{\pi_{k,k'}\}} \prob_{\mu,\pi^*} \big( \hat\pi\neq\pi^*\big) - \mu(\mathcal E\setminus\bar\T_\k) \\
&\ge \frac{\sqrt{15}}{\sqrt{15}+1} \Big( \frac{3}{4} - \frac{1}{2\sqrt{\log15}} \Big) - \frac{m(m-1)}2 e^{-d/8},
\end{align}
where in the last inequality we applied the result of \Cref{lower-bound} in conjunction with the monotonicity of function $m \mapsto \frac{\sqrt{m}}{1+\sqrt{m}}(3/4 - (2\sqrt{\log(m)})^{-1})$. Recall that $m > n \ge 5$ and $d \ge 16 \log(nm)$ yielding $\inf_{\hat\pi} \prob_{\btdiese, \bsigmadiese, \pi^*}(\hat\pi \neq \pi^*) > 0.338$, concluding the proof.

\subsection{Proof of \Cref{thm:mild}}
    To ease notation, we write $\hat\pi$ instead
    of $\hat\pi^{\textup{LSL}}_{n,m}$, and, without
    loss of generality, we assume that $\pi^*(i) 
    = i$ for $i\in [n]$. We wish to 
    prove that on an event of probability $\ge 1
    -\alpha$, for every injective mapping $\pi
    : [n]\to [m]$, we have $\psi(\pi^*)\le\psi 
    (\pi)$, where
    \begin{align}
        \psi(\pi) = \sum_{i=1}^n \log \|X_i - 
        \Xdiese_{\pi(i)}\|^2.
    \end{align}
    Since the logarithm is an increasing function, 
    this is equivalent to showing that
    \begin{align}
        \prod_{i=1}^n  \|X_i-\Xdiese_{\pi^*(i)}\|^2 
        & < \prod_{i=1}^n \|X_i-\Xdiese_{\pi(i)}\|^2,
        \quad\text{for every $\pi\neq\pi^*$},
    \end{align}
    which, in turn, is the same as 
    \begin{align}
        \prod_{i=1}^n  \frac{\|X_i-\Xdiese_{
        \pi^*(i)}\|^2}{\|X_i-\Xdiese_{\pi(i)} \|^2} 
        <1, \quad\text{for every $\pi\neq\pi^*$}.
    \end{align}
    In view of \eqref{upbX} and \eqref{lowbX}, 
    we have
    \begin{align}
        \prod_{i=1}^n  \frac{ \|X_i - \Xdiese_{ 
        \pi^*(i)}\|^2 }{ \|X_i - \Xdiese_{\pi(i)}
        \|^2} & \le 
        \prod_{\substack{i\in [n]\\ \pi(i) \neq 
        \pi^*(i)}}  \frac{ 2\sigma_i^2 (d +
        \sqrt{d}\,\zeta_2)}{ \sigma_{i,\pi(i)}^2 
        \big ( \kappa_{i, \pi(i)}^2 + d - \sqrt{d} 
        \,\zeta_2 - 2\zeta_1 \kappa_{i,\pi(i)} 
        \big)_+}\\
        &\le \prod_{\substack{i\in [n]\\ \pi(i) \neq 
        \pi^*(i)}}  \frac{ 4\sigma_i^2 (d +
        \sqrt{d}\,\zeta_2) }{ \sigma_{i,\pi(i)}^2 
        \big(\kappa_{i,\pi(i)}^2 + 2d - 2\sqrt{d}\,
        \zeta_2 \big)_+},\qquad \text{if }\quad 
        \zeta_1 \le (1/4)\bar\kappa.\label{prod0}
    \end{align}
    Let us define the sets $I_1 = \big\{i\in 
    [n]: \pi(i) \in\Im(\pi^*)\setminus\{\pi^*(i)
    \}\big\}$ and $I_2 = \{i\in [n]: \pi(i) 
    \not\in\Im(\pi^*)\}$. Clearly, using the 
    inequality $\sigma_{i,j}^2 \ge 2\sigma_i 
    \sdiese_j $, we get
    \begin{align}\label{prod1}
        \prod_{\substack{i\in [n]\\ \pi(i) \neq 
        \pi^*(i)}}  \frac{ 2\sigma_i^2}{ 
        \sigma_{i,\pi(i)}^2}
        &\le \prod_{\substack{i\in [n]\\ \pi(i) 
        \neq \pi^*(i)}}  \frac{ \sigma_i^2}{ 
        \sigma_i\sdiese_{\pi(i)}}
        =  \frac{ \prod_{i\in I_1\cup I_2} 
        \sigma_i}{\prod_{i\in I_1}  
        \sdiese_{\pi(i)} \prod_{i\in I_2}  
        \sdiese_{\pi(i)}}.
    \end{align}
    For every $i\in I_1$, there is $j\in [n]$ 
    such that $\pi(i) = \pi^*(j)$; this $j$ is 
    given by $j = (\pi^*)^{-1}(i)$. For such a 
    pair $(i,j)$, in view of \eqref{pistar}, 
    we have $\sdiese_{\pi(i)} = 
    \sdiese_{\pi^*(j)} = \sigma_j$. Note that 
    by construction of $I_1$, $(\pi^*)^{-1}(I_1) 
    \subset I_1\cup I_2$. This implies that 
    \begin{align}\label{prod2}
        \prod_{i\in I_1}  \sdiese_{\pi(i)} = 
        \prod_{j\in (\pi^*)^{-1}(I_1)}  \sigma_j
        = \frac{\prod_{j\in I_1\cup I_2} \sigma_j
        }{\prod_{j\in (I_1\cup I_2)\setminus
        (\pi^*)^{-1} (I_1)}\sigma_j}.
    \end{align}
    Note also that the cardinality of the set 
    $J_1 = (\pi^*)^{-1} (I_1)$ is equal to 
    the cardinality of $I_1$, which implies that
    $|(I_1\cup I_2)\setminus J_1| = |I_2|$. 
    Combining \eqref{prod1}, \eqref{prod2}, 
    and the last equality of cardinalities, 
    we get 
    \begin{align}
        \prod_{\substack{i\in [n]\\ \pi(i) \neq 
        \pi^*(i)}}  \frac{ 2\sigma_i^2}{ \sigma_{
        i,\pi(i)}^2}
        &\le\frac{\prod_{j\in (I_1\cup I_2)
        \setminus J_1}\sigma_j}{\prod_{i\in I_2}
        \sdiese_{\pi(i)}}\le r_\sigma^{|I_2|}.
        \label{prod3}
    \end{align}
    Using the same notation $I_1$ and $I_2$, 
    we can check that 
    \begin{align}
        \kappa_{i,\pi(i) }
        \ge 
        \begin{cases}
        \kinin,& i\in I_1,\\
        \kinout,& i\in I_2.
        \end{cases}
    \end{align}
    Injecting this inequality into \eqref{prod0}, 
    and using \eqref{prod3}, we get 
    \begin{align}
        \prod_{i=1}^n  \frac{ \|X_i - \Xdiese_{ 
        \pi^*(i)}\|^2 }{ \|X_i - \Xdiese_{\pi(i)}
        \|^2} & \le \frac{ r_{\sigma}^{I_2} \{2(d 
        + \sqrt{d}\, \zeta_2)\}^{|I_1| + |I_2|} }{
        (\kinin^2 + 2d -
        2\sqrt{d}\,\zeta_2)_+^{I_1} (\bar\kappa_{
        \textup{in-out}}^2 + 2d -
        2\sqrt{d}\,\zeta_2)_+^{I_2}}.
    \end{align}
    Recall that this inequality is true on the event
    $\zeta_1\le \bar\kappa/4$. It follows from last
    display that as soon as
    \begin{align}\label{ineq:kappa}
    \begin{cases}
        \zeta_1&\le \bar\kappa/4\\
        4\sqrt{d}\,\zeta_2 &<  \bar\kappa_{
        \textup{in-in}}^2\\
        2d(r_\sigma-1) + 4r_\sigma\sqrt{d}\,\zeta_2
        &\le \kinout^2
    \end{cases}
    \end{align}
    we have 
    \begin{align}
        \prod_{i=1}^n  \frac{ \|X_i - \Xdiese_{ 
        \pi^*(i)}\|^2 }{ \|X_i - \Xdiese_{\pi(i)}
        \|^2} <1
    \end{align}
    for every $\pi$. It remains to show that,
    under the conditions of \Cref{thm:mild}, the
    event in \eqref{ineq:kappa} has a probability
    at least $1-\alpha$. This will be done by 
    using tail bounds for Gaussian and $\chi$-squared
    distributions, combined with the union bound. 
    
    On the one hand, using the well-known tail 
    bound for the standard Gaussian distribution 
    and the union bound, we get
    \begin{align}\label{z1-bound:2}
    \prob\bigg(\z_1\ge \sqrt{2\log\Big(\frac{4nm}{
    \alpha}\Big)}\bigg)
    &\le \sum_{i\not =j} \prob\bigg(|\z_{i,j}|\ge
    \sqrt{2\log\Big(\frac{4nm}{\alpha}\Big)}\bigg)
    \le \alpha/2.
    \end{align}
    On the other hand, \Cref{concentration} and the
    union bound entail
    \begin{align}\label{z2-bound:2}
    \prob\bigg(\z_2\ge 2\sqrt{\log({4nm}/{\alpha})} 
    + \frac{2\log(4nm/\alpha)}{\sqrt{d}}\bigg)
    &\le \alpha/2.
    \end{align}
    Therefore, if 
    \begin{align}
        \begin{cases}
            \bar\kappa &\ge\ 4 \sqrt{2\log(4nm/
            \alpha)}\\
            \kinin^2 &\ge\  
            8\sqrt{d\log(4nm/\alpha)} +
            8{\log(4nm/\alpha)}\\
            \kinout^2 &\ge\  
            2d(r_\sigma-1) + 8r_\sigma
            \sqrt{d\log(4nm/\alpha)}\,+
            8r_\sigma{\log(4nm/\alpha)}
        \end{cases}
    \end{align}
    then, on an event of probability $\ge 1-\alpha$, 
    all the inequalities in \eqref{ineq:kappa} 
    hold true.  This completes the proof of the
    theorem.

\begin{acks}[Acknowledgments]
	This work was supported by the grant
	Investissements d'Avenir (ANR-11-IDEX-0003/Labex Ecodec/ANR-11-LABX-0047) and by the FAST Advance grant.
\end{acks}
\bibliography{literature}

\begin{thebibliography}{}

\bibitem[\protect\astroncite{Azaïs and de~Castro}{2020}]{azais2020multiple}
Azaïs, J.~M. and de~Castro, Y. (2020).
\newblock Multiple testing and variable selection along least angle
  regression's path.

\bibitem[\protect\astroncite{Bai et~al.}{2020}]{bai2020d3feat}
Bai, X., Luo, Z., Zhou, L., Fu, H., Quan, L., and Tai, C.-L. (2020).
\newblock D3feat: Joint learning of dense detection and description of 3d local
  features.
\newblock In {\em Proceedings of the IEEE/CVF Conference on Computer Vision and
  Pattern Recognition (CVPR)}.

\bibitem[\protect\astroncite{{Blanchard} et~al.}{2018}]{Balchard}
{Blanchard}, G., {Carpentier}, A., and {Gutzeit}, M. (2018).
\newblock {Minimax Euclidean separation rates for testing convex hypotheses in
  \(\mathbb{R}^{d}\)}.
\newblock {\em {Electron. J. Stat.}}, 12(2):3713--3735.

\bibitem[\protect\astroncite{{Burnashev}}{1979}]{Burnashev}
{Burnashev}, M.~V. (1979).
\newblock {On the minimax detection of an inaccurately known signal in a white
  Gaussian noise background}.
\newblock {\em {Theory Probab. Appl.}}, 24:107--119.

\bibitem[\protect\astroncite{Calonder et~al.}{2010}]{brief2010}
Calonder, M., Lepetit, V., Strecha, C., and Fua, P. (2010).
\newblock Brief: Binary robust independent elementary features.
\newblock In Daniilidis, K., Maragos, P., and Paragios, N., editors, {\em
  Computer Vision -- ECCV 2010}, pages 778--792, Berlin, Heidelberg. Springer
  Berlin Heidelberg.

\bibitem[\protect\astroncite{Carpentier et~al.}{2019}]{carpentier2019minimax}
Carpentier, A., Collier, O., Comminges, L., Tsybakov, A.~B., and Wang, Y.
  (2019).
\newblock Minimax rate of testing in sparse linear regression.
\newblock {\em Automation and Remote Control}, 80(10):1817--1834.

\bibitem[\protect\astroncite{Chen et~al.}{2010}]{piifd2010}
Chen, J., Tian, J., Lee, N., Zheng, J., Smith, R.~T., and Laine, A.~F. (2010).
\newblock A partial intensity invariant feature descriptor for multimodal
  retinal image registration.
\newblock {\em IEEE Transactions on Biomedical Engineering}, 57(7):1707--1718.

\bibitem[\protect\astroncite{Collier}{2012}]{Collier12a}
Collier, O. (2012).
\newblock Minimax hypothesis testing for curve registration.
\newblock In {\em {AISTATS} 2012}, volume~22 of {\em {JMLR} Proceedings}, pages
  236--245.

\bibitem[\protect\astroncite{Collier and Dalalyan}{2013}]{jmlr_CD13}
Collier, O. and Dalalyan, A.~S. (2013).
\newblock Permutation estimation and minimax rates of identifiability.
\newblock {\em Journal of Machine Learning Research}, W {{\&}} CP 31 (AI-STATS
  2013):10--19.

\bibitem[\protect\astroncite{Collier and Dalalyan}{2016}]{collier2016minimax}
Collier, O. and Dalalyan, A.~S. (2016).
\newblock Minimax rates in permutation estimation for feature matching.
\newblock {\em The Journal of Machine Learning Research}, 17(1):162--192.

\bibitem[\protect\astroncite{{Comminges} and {Dalalyan}}{2012}]{ComDal1}
{Comminges}, L. and {Dalalyan}, A.~S. (2012).
\newblock {Tight conditions for consistency of variable selection in the
  context of high dimensionality}.
\newblock {\em {Ann. Stat.}}, 40(5):2667--2696.

\bibitem[\protect\astroncite{{Comminges} and {Dalalyan}}{2013}]{ComDal2}
{Comminges}, L. and {Dalalyan}, A.~S. (2013).
\newblock {Minimax testing of a composite null hypothesis defined via a
  quadratic functional in the model of regression}.
\newblock {\em {Electron. J. Stat.}}, 7:146--190.

\bibitem[\protect\astroncite{Duff and Koster}{2001}]{Duff}
Duff, I.~S. and Koster, J. (2001).
\newblock On algorithms for permuting large entries to the diagonal of a sparse
  matrix.
\newblock {\em SIAM Journal on Matrix Analysis and Applications},
  22(4):973--996.

\bibitem[\protect\astroncite{Ermakov}{1990}]{Ermakov}
Ermakov, M.~S. (1990).
\newblock Minimax detection of a signal in {G}aussian white noise.
\newblock {\em Teor. Veroyatnost. i Primenen.}, 35(4):704--715.

\bibitem[\protect\astroncite{Flammarion et~al.}{2019}]{flammarion2019optimal}
Flammarion, N., Mao, C., and Rigollet, P. (2019).
\newblock Optimal rates of statistical seriation.
\newblock {\em Bernoulli}, 25(1):623--653.

\bibitem[\protect\astroncite{Gao and Zhang}{2019}]{gao2019iterative}
Gao, C. and Zhang, A.~Y. (2019).
\newblock Iterative algorithm for discrete structure recovery.
\newblock {\em arXiv preprint arXiv:1911.01018}.

\bibitem[\protect\astroncite{Harwood and Drummond}{2016}]{harwood2016}
Harwood, B. and Drummond, T. (2016).
\newblock Fanng: Fast approximate nearest neighbour graphs.
\newblock In {\em 2016 IEEE Conference on Computer Vision and Pattern
  Recognition (CVPR)}, pages 5713--5722.

\bibitem[\protect\astroncite{{Ingster}}{1982}]{Ingster82}
{Ingster}, Y.~I. (1982).
\newblock {Minimax nonparametric detection of signals in white Gaussian noise}.
\newblock {\em {Probl. Inf. Transm.}}, 18:130--140.

\bibitem[\protect\astroncite{Ingster and Suslina}{2003}]{Ingster}
Ingster, Y.~I. and Suslina, I.~A. (2003).
\newblock {\em Nonparametric goodness-of-fit testing under {G}aussian models},
  volume 169 of {\em Lecture Notes in Statistics}.
\newblock Springer-Verlag, New York.

\bibitem[\protect\astroncite{Itseez}{2015}]{opencv}
Itseez (2015).
\newblock Open source computer vision library.
\newblock \url{https://github.com/itseez/opencv}.

\bibitem[\protect\astroncite{Jiang et~al.}{2016}]{jiang2016}
Jiang, Z., Xie, L., Deng, X., Xu, W., and Wang, J. (2016).
\newblock Fast nearest neighbor search in the hamming space.
\newblock In Tian, Q., Sebe, N., Qi, G.-J., Huet, B., Hong, R., and Liu, X.,
  editors, {\em MultiMedia Modeling}, pages 325--336, Cham. Springer
  International Publishing.

\bibitem[\protect\astroncite{Jin et~al.}{2020}]{Jin2020}
Jin, Y., Mishkin, D., Mishchuk, A., Matas, J., Fua, P., Yi, K.~M., and Trulls,
  E. (2020).
\newblock {Image Matching across Wide Baselines: From Paper to Practice}.
\newblock {\em International Journal of Computer Vision}.

\bibitem[\protect\astroncite{{Juditsky} and {Nemirovski}}{2020}]{JudNem}
{Juditsky}, A. and {Nemirovski}, A. (2020).
\newblock {\em {Statistical inference via convex optimization}}.
\newblock Princeton, NJ: Princeton University Press.

\bibitem[\protect\astroncite{Kuhn}{1955}]{Kuhn}
Kuhn, H.~W. (1955).
\newblock {The Hungarian Method for the Assignment Problem}.
\newblock {\em Naval Research Logistics Quarterly}, 2(1--2):83--97.

\bibitem[\protect\astroncite{Kuhn}{2010}]{Kuhn2010}
Kuhn, H.~W. (2010).
\newblock {\em The Hungarian Method for the Assignment Problem}, pages 29--47.
\newblock Springer Berlin Heidelberg, Berlin, Heidelberg.

\bibitem[\protect\astroncite{Laurent and Massart}{2000}]{LaurentMassart2000}
Laurent, B. and Massart, P. (2000).
\newblock Adaptive estimation of a quadratic functional by model selection.
\newblock {\em Ann. Statist.}, 28(5):1302--1338.

\bibitem[\protect\astroncite{Lepski and
  Tsybakov}{2000}]{lepski2000asymptotically}
Lepski, O.~V. and Tsybakov, A.~B. (2000).
\newblock Asymptotically exact nonparametric hypothesis testing in sup-norm and
  at a fixed point.
\newblock {\em Probability Theory and Related Fields}, 117(1):17--48.

\bibitem[\protect\astroncite{Lowe}{2004}]{lowe2004distinctive}
Lowe, D.~G. (2004).
\newblock Distinctive image features from scale-invariant keypoints.
\newblock {\em International journal of computer vision}, 60(2):91--110.

\bibitem[\protect\astroncite{Ma et~al.}{2021}]{ma2021}
Ma, J., Jiang, X., Fan, A., Jiang, J., and Yan, J. (2021).
\newblock Image matching from handcrafted to deep features: A survey.
\newblock {\em International Journal of Computer Vision}, 129.

\bibitem[\protect\astroncite{Ma et~al.}{2020}]{ma2020optimal}
Ma, R., Tony~Cai, T., and Li, H. (2020).
\newblock Optimal permutation recovery in permuted monotone matrix model.
\newblock {\em Journal of the American Statistical Association}, pages 1--15.

\bibitem[\protect\astroncite{Malkov and Yashunin}{2020}]{malkov2020}
Malkov, Y.~A. and Yashunin, D.~A. (2020).
\newblock Efficient and robust approximate nearest neighbor search using
  hierarchical navigable small world graphs.
\newblock {\em IEEE Transactions on Pattern Analysis and Machine Intelligence},
  42(4):824--836.

\bibitem[\protect\astroncite{Mao et~al.}{2020}]{Mao18}
Mao, C., Pananjady, A., and Wainwright, M.~J. (2020).
\newblock Towards optimal estimation of bivariate isotonic matrices with
  unknown permutations.
\newblock {\em Annals of Statistics}, 48(6):3183--3205.

\bibitem[\protect\astroncite{Mao et~al.}{2018}]{mao2018minimax}
Mao, C., Weed, J., and Rigollet, P. (2018).
\newblock Minimax rates and efficient algorithms for noisy sorting.
\newblock In {\em Algorithmic Learning Theory}, pages 821--847. PMLR.

\bibitem[\protect\astroncite{Munkres}{1957}]{Munkres}
Munkres, J. (1957).
\newblock Algorithms for the assignment and transportation problems.
\newblock {\em Journal of the Society for Industrial and Applied Mathematics},
  5(1):32--38.

\bibitem[\protect\astroncite{{Ndaoud} and {Tsybakov}}{2020}]{Ndaoud}
{Ndaoud}, M. and {Tsybakov}, A.~B. (2020).
\newblock {Optimal variable selection and adaptive noisy compressed sensing}.
\newblock {\em {IEEE Trans. Inf. Theory}}, 66(4):2517--2532.

\bibitem[\protect\astroncite{Pananjady and
  Samworth}{2020}]{pananjady2020isotonic}
Pananjady, A. and Samworth, R.~J. (2020).
\newblock Isotonic regression with unknown permutations: Statistics,
  computation, and adaptation.
\newblock {\em arXiv preprint arXiv:2009.02609}.

\bibitem[\protect\astroncite{Pananjady et~al.}{2017}]{pananjady2017linear}
Pananjady, A., Wainwright, M.~J., and Courtade, T.~A. (2017).
\newblock Linear regression with shuffled data: Statistical and computational
  limits of permutation recovery.
\newblock {\em IEEE Transactions on Information Theory}, 64(5):3286--3300.

\bibitem[\protect\astroncite{Ramdas et~al.}{2016}]{RamdasISW16}
Ramdas, A., Isenberg, D., Singh, A., and Wasserman, L.~A. (2016).
\newblock Minimax lower bounds for linear independence testing.
\newblock In {\em {IEEE} International Symposium on Information Theory, {ISIT}
  2016, Barcelona, Spain, July 10-15, 2016}, pages 965--969. {IEEE}.

\bibitem[\protect\astroncite{Rublee et~al.}{2011}]{orb2011}
Rublee, E., Rabaud, V., Konolige, K., and Bradski, G. (2011).
\newblock Orb: An efficient alternative to sift or surf.
\newblock In {\em 2011 International Conference on Computer Vision}, pages
  2564--2571.

\bibitem[\protect\astroncite{Shah et~al.}{2021}]{ShahBW16a}
Shah, N.~B., Balakrishnan, S., and Wainwright, M.~J. (2021).
\newblock A permutation-based model for crowd labeling: Optimal estimation and
  robustness.
\newblock {\em IEEE Transactions on Information Theory}, 67(6):4162--4184.

\bibitem[\protect\astroncite{Slawski and Ben-David}{2019}]{slawski2019linear}
Slawski, M. and Ben-David, E. (2019).
\newblock Linear regression with sparsely permuted data.
\newblock {\em Electronic Journal of Statistics}, 13(1):1--36.

\bibitem[\protect\astroncite{Tian et~al.}{2020}]{Tian_2020_ACCV}
Tian, Y., Balntas, V., Ng, T., Barroso-Laguna, A., Demiris, Y., and
  Mikolajczyk, K. (2020).
\newblock D2d: Keypoint extraction with describe to detect approach.
\newblock In {\em Proceedings of the Asian Conference on Computer Vision
  (ACCV)}.

\bibitem[\protect\astroncite{Tsybakov}{2009}]{Tsybakov2009}
Tsybakov, A.~B. (2009).
\newblock {\em Introduction to nonparametric estimation}.
\newblock Springer Series in Statistics. Springer, New York.

\bibitem[\protect\astroncite{Virtanen et~al.}{2020}]{2020SciPy-NMeth}
Virtanen, P., Gommers, R., Oliphant, T.~E., Haberland, M., Reddy, T.,
  Cournapeau, D., Burovski, E., Peterson, P., Weckesser, W., Bright, J., {van
  der Walt}, S.~J., Brett, M., Wilson, J., Millman, K.~J., Mayorov, N., Nelson,
  A. R.~J., Jones, E., Kern, R., Larson, E., Carey, C.~J., Polat, {\.I}., Feng,
  Y., Moore, E.~W., {VanderPlas}, J., Laxalde, D., Perktold, J., Cimrman, R.,
  Henriksen, I., Quintero, E.~A., Harris, C.~R., Archibald, A.~M., Ribeiro,
  A.~H., Pedregosa, F., {van Mulbregt}, P., and {SciPy 1.0 Contributors}
  (2020).
\newblock {{SciPy} 1.0: Fundamental Algorithms for Scientific Computing in
  Python}.
\newblock {\em Nature Methods}, 17:261--272.

\bibitem[\protect\astroncite{Wang et~al.}{2018}]{wang2018}
Wang, K., Zhu, N., Cheng, Y., Li, R., Zhou, T., and Long, X. (2018).
\newblock Fast feature matching based on r -nearest k -means searching.
\newblock {\em CAAI Transactions on Intelligence Technology}, 3(4):198--207.

\bibitem[\protect\astroncite{Wang}{2011}]{wang2011}
Wang, X. (2011).
\newblock A fast exact k-nearest neighbors algorithm for high dimensional
  search using k-means clustering and triangle inequality.
\newblock In {\em The 2011 International Joint Conference on Neural Networks},
  pages 1293--1299.

\bibitem[\protect\astroncite{Wei and Wainwright}{2020}]{Yuting}
Wei, Y. and Wainwright, M.~J. (2020).
\newblock The local geometry of testing in ellipses: Tight control via
  localized kolmogorov widths.
\newblock {\em IEEE Transactions on Information Theory}, 66(8):5110--5129.

\bibitem[\protect\astroncite{Wei et~al.}{2019a}]{Yuting2}
Wei, Y., Wainwright, M.~J., and Guntuboyina, A. (2019a).
\newblock {The geometry of hypothesis testing over convex cones: Generalized
  likelihood ratio tests and minimax radii}.
\newblock {\em The Annals of Statistics}, 47(2):994 -- 1024.

\bibitem[\protect\astroncite{Wei et~al.}{2019b}]{Wei}
Wei, Y., Wainwright, M.~J., and Guntuboyina, A. (2019b).
\newblock {The geometry of hypothesis testing over convex cones: Generalized
  likelihood ratio tests and minimax radii}.
\newblock {\em The Annals of Statistics}, 47(2):994 -- 1024.

\bibitem[\protect\astroncite{Wolfer and Kontorovich}{2020}]{Wolfer}
Wolfer, G. and Kontorovich, A. (2020).
\newblock Minimax testing of identity to a reference ergodic markov chain.
\newblock In {\em {AISTATS} 2020}, volume 108 of {\em Proceedings of Machine
  Learning Research}, pages 191--201.

\bibitem[\protect\astroncite{Xing et~al.}{2020}]{Xing20}
Xing, X., Liu, M., Ma, P., and Zhong, W. (2020).
\newblock Minimax nonparametric parallelism test.
\newblock {\em Journal of Machine Learning Research}, 21(94):1--47.

\end{thebibliography}

\end{document}